\documentclass[12pt]{article}
\usepackage{xcolor}
\usepackage{amsfonts}
\usepackage{amssymb}
\usepackage{CJK}
\usepackage{cite}
\usepackage{setspace}

\usepackage[utf8]{inputenc}
\usepackage{amsthm}
\usepackage{amsmath,bm}
\usepackage{graphicx}
\usepackage[citecolor=red]{hyperref}
\usepackage{tikz}
\usepackage{units}
\usepackage{mathrsfs}
\usepackage{stmaryrd}
\usepackage{mathabx}
\usetikzlibrary{decorations.pathmorphing}
\usepackage{indentfirst}
\usepackage{mdframed}
\allowdisplaybreaks
\numberwithin{equation}{section}

\makeatletter
\@namedef{subjclassname@2020}{%
	\textup{2020} Mathematics Subject Classification}
\makeatother

\newcommand{\beq}{\begin{equation}}
\newcommand{\eeq}{\end{equation}}

\newtheorem{thm}{Theorem}[section]

\newtheorem{lem}[thm]{Lemma}
\newtheorem{pro}[thm]{Proposition}
\theoremstyle{definition}

\usepackage[bottom]{footmisc}

\title{Stability for the Boussinesq Equations with Horizontal Dissipation near the Hydrostatic Balance on $\mathbb{R}^2$}

\author{
	Jiahong Wu\footnote{Department of Mathematics, University of Notre Dame, Notre Dame, IN 46556, USA. Email: jwu29@nd.edu}
	\and
	Mengxin Yan\footnote{School of Mathematics, Shandong University, Jinan, Shandong 250000, P. R. China. Email: 3020135975@qq.com}
	\and
	Ning Zhu\footnote{School of Mathematics, Shandong University, Jinan, Shandong 250000, P. R. China. Email: ning.zhu@sdu.edu.cn}
}
\date{}

\begin{document}
	
	\maketitle
	
	\begin{abstract}
		The hydrostatic balance is a fundamental equilibrium state in stratified fluids and plays a central role in geophysical fluid dynamics. Understanding its stability under incomplete dissipation is a longstanding challenge, since anisotropic diffusion alone is generally insufficient to control the nonlinear evolution and no robust stabilizing mechanism is known for the corresponding anisotropically dissipative Navier--Stokes equations. In this paper, we investigate the two-dimensional Boussinesq equations on $\mathbb{R}^2$ with only horizontal dissipation near the hydrostatic equilibrium $(U,\Theta)=(0,x_2)$. We show that the velocity--temperature coupling generates internal gravity waves whose dispersive decay, together with the horizontal dissipation, provides an effective stabilizing mechanism that compensates for the complete absence of vertical dissipation. This identifies a mechanism by which dispersive wave propagation restores stability in an incompletely dissipative fluid system. For sufficiently small initial perturbations in $H^k(\mathbb{R}^2)\cap W^{3,1}(\mathbb{R}^2)$ with $k\ge14$, we establish the global existence and uniqueness of classical solutions together with explicit anisotropic, componentwise large-time decay rates for the velocity and temperature, including faster decay of the vertical velocity.
	\end{abstract}
	
	

\begin{spacing}{1.2}
	\section{Introduction}
	\label{intr}
	
A fundamental problem in fluid dynamics is to understand whether small perturbations of an equilibrium remain stable and eventually decay as time evolves. For many important fluid models, such a stabilizing mechanism is either absent or still poorly understood. Even with small initial data, the three-dimensional (3D) incompressible Euler equations, for example, may develop finite-time singularities, as demonstrated in recent breakthrough works by Elgindi \cite{Elgindi}, Chen--Hou \cite{ChenHou1,ChenHou2}, C\'ordoba \cite{CMZ}, and Shao--Wei--Zhang \cite{SWZ}. Even in two dimensions, where global regularity is known, solutions of the incompressible Euler equations may experience extremely rapid growth of Sobolev norms, with double exponential growth in time being possible (see, e.g., \cite{KS,Zlatos,CJ,JYZ}). 
	
\vskip .1in 
On the other hand, for the incompressible Navier--Stokes equations with the full Laplacian dissipation, every finite-energy weak solution in the whole space decays algebraically to zero as $t\to\infty$ (see Schonbek \cite{Schonbek}). However, once only partial dissipation is present, the situation becomes much less clear. Consider, for example, the two-dimensional Navier--Stokes equations with only horizontal dissipation,
\begin{equation}\label{equation01}
	\left\{
		\begin{aligned}
			&\partial_tu-\partial_1^2u+u\cdot\nabla u+\nabla p=0,
			\qquad x\in\mathbb R^2,\ t>0,\\
			&\nabla\cdot u=0,\\
			&u(x,0)=u_0(x).
		\end{aligned}
		\right.
\end{equation}	
For sufficiently smooth initial data, for instance $u_0\in H^2(\mathbb R^2)$, the above system admits a unique global solution. The real issue is not global existence but rather its long-time behavior: does every solution remain uniformly bounded and converge to zero, or can higher Sobolev norms grow? Standard energy estimates yield at best a doubly exponential upper bound on the $H^2$ norm, providing essentially no information on the actual asymptotic behavior. In a recent work \cite{WangWuZhu}, Wang, Wu, and Zhu proved that replacing $\partial_1^2$ by the fractional operator $(-\partial_1^2)^\alpha$ with $0\le\alpha<1$ leads to global decay of solutions. However, the physically important endpoint case $\alpha=1$ remains completely open.
	
\vskip .1in 
A remarkable phenomenon emerges when incompletely dissipative fluid equations are coupled with additional physical effects such as buoyancy, rotation, or a background magnetic field. Unlike the underlying fluid equations, the coupled systems possess an intrinsic hyperbolic wave structure. The propagation of these waves gives rise to dispersive effects, which spread the energy over space and generate temporal decay of the solution. Although this decay is fundamentally different from the exponential damping produced by viscosity, it nevertheless provides an effective stabilizing mechanism. When combined with even weak or anisotropic dissipation, the dispersive decay can produce additional smoothing and time integrability that are unavailable from the dissipative terms alone. Understanding how such weak wave-induced stabilization compensates for the lack of full dissipation has become a central theme in the mathematical analysis of geophysical fluids and magnetohydrodynamics.

\vskip .1in 
The Boussinesq equations constitute one of the simplest and most important fluid models exhibiting this stabilization mechanism. The coupling between the velocity and the temperature naturally generates internal gravity waves, whose dispersive decay can substantially alter the long-time behavior of the fluid. Motivated by this phenomenon, we consider the two-dimensional Boussinesq equations with only horizontal dissipation,
	\begin{align}\label{equation0}
	\left \{
	\begin{aligned} 
		& \partial_t U  -  \partial_1^2U +U \cdot \nabla U  + \nabla P     =   \Theta e_2 ,&   x \in \mathbb{R}^2, \ t > 0,\\
		& \partial_t \Theta  -  \partial_1^2 \Theta  + U \cdot \nabla \Theta = 0     ,\\
		& \nabla \cdot U =0    ,  \\
		& U(x,0)= U_0(x), \  \Theta(x,0) = \Theta_0(x)   , 
	\end{aligned}
	\right.
\end{align}
where $U(x,t) = (U_1(x,t), U_2(x,t))$ denotes the fluid velocity, $\Theta(x,t)$ the temperature, $P(x,t)$ the scalar pressure. $\partial_1$ is an abbreviation of $\partial_{x_1}$ and $e_2=(0,1)$, the forcing term represents the buoyancy force of the fluid. 

\vskip .1in
The Boussinesq equations are fundamental models for buoyancy-driven incompressible flows and arise naturally in atmospheric dynamics, oceanography, and Rayleigh--Bénard convection. In these models, the buoyancy force generated by temperature variations interacts with the incompressible velocity field to produce internal gravity waves, whose dispersive properties play a fundamental role in the large-time dynamics. Systems with only partial or anisotropic dissipation also arise in several physically relevant situations, with Prandtl's boundary layer equations being one of the most notable examples.

\vskip .1in
Among the equilibrium states of the Boussinesq equations, the hydrostatic balance plays a distinguished role in geophysical fluid dynamics. It describes a fluid at rest in which the vertical pressure gradient exactly balances the buoyancy force. Such a balance provides the leading-order approximation for many large-scale atmospheric and oceanic flows and serves as a fundamental reference state in the study of stratified fluids.
More precisely,
\[
U_{he}\triangleq(0,0),\qquad
\Theta_{he}\triangleq x_2,\qquad
P_{he}\triangleq\frac12x_2^2
\]
forms a steady solution of \eqref{equation0}. The main objective of this paper is to investigate the nonlinear stability and large-time behavior of solutions that are initially close to this hydrostatic equilibrium.

\vskip .1in 
To this end, we introduce the perturbation variables
\[
u=U-U_{he},\qquad
\theta=\Theta-\Theta_{he},\qquad
p=P-P_{he}.
\]
Substituting these expressions into \eqref{equation0}, one easily verifies that the perturbation $(u,\theta,p)$ satisfies
\begin{align}\label{equation1}
	\left \{
	\begin{aligned} 
		& \partial_t u  -  \partial_1^2u +u \cdot \nabla u  + \nabla p     =    \theta e_2, &   x \in \mathbb{R}^2, \quad  t > 0,\\
		& \partial_t \theta  - \partial_1^2 \theta + u_2  + u \cdot \nabla \theta = 0     ,\\
		& \nabla \cdot u =0    ,  \\
		& u(0,x)= u_0(x), \  \theta(0,x) = \theta_0(x)   .
	\end{aligned}
	\right.
\end{align}
The perturbation system \eqref{equation1} exhibits a remarkable feature absent from the incompressible Navier--Stokes equations. Owing to the coupling between the buoyancy force $\theta e_2$ and the stratification term $u_2$, the linearized system possesses an intrinsic internal gravity wave structure. As will be shown below, the dispersive decay associated with these waves provides an effective stabilizing mechanism that compensates for the lack of full dissipation.

\vskip .1in 
The stability problem for the perturbation system \eqref{equation1} cannot be resolved by standard energy methods. As discussed above, the horizontal dissipation fails to provide sufficient smoothing for the velocity field and the classical energy method  yields at best a rapidly growing upper bound for higher Sobolev norms. Treating the velocity equation separately therefore fails to capture the mechanism responsible for stabilization. The essential difficulty stems from the lack of vertical dissipation. At the level of the vorticity ($\omega=\nabla\times u$) equation, the nonlinear terms contain products involving vertical derivatives that cannot be absorbed by the available horizontal dissipation. In particular, the most delicate nonlinear interactions 
$$
\int_{\mathbb{R}^2} \partial_2 u_2 (\partial_2 \omega)^2 dx
$$
generate coefficients that, under the decay expected from the one-dimensional heat operator, behave only like $(1+t)^{-1}$, which is critical and fails to provide the time integrability needed to close the energy estimates. Consequently, the standard energy method alone cannot establish uniform-in-time bounds, let alone decay of the solution.

\vskip .1in 
The key idea of this paper is to regard \eqref{equation1} as a genuinely coupled system rather than as the Navier--Stokes equations with a forcing term. The coupling between the momentum and temperature equations gives rise to a wave structure that is invisible from the individual equations. More precisely, the buoyancy force in the momentum equation and the stratification term $u_2$ in the temperature equation combine to generate internal gravity waves. These waves exhibit dispersive decay, which weakens the nonlinear interactions and provides an effective damping mechanism. Although this wave-induced stabilization is substantially weaker than the classical parabolic smoothing generated by the full Laplacian, it supplies precisely the additional decay needed to compensate for the missing dissipation.

\vskip .1in 
The main challenge is to quantify this dispersive stabilization at the nonlinear level and to combine it with the available horizontal dissipation in a manner that yields uniform energy bounds and ultimately global nonlinear stability. Our main result shows that this program can indeed be carried out. More precisely, we establish the global stability of the hydrostatic equilibrium together with explicit anisotropic decay estimates for the perturbation.
\begin{thm}\label{main2}
Let $k\ge14$. Consider \eqref{equation1} with the initial data $u_0, \theta_0 \in H^{k}(\mathbb{R}^2) \bigcap  W^{3,1}(\mathbb{R}^2)$  and   
	$\nabla \cdot u_0 = 0$. Then there exists $\varepsilon>0$ such that whenever
	\begin{equation}\label{Initial data condition}
	\|u_0\|_{H^k}
	+\|\theta_0\|_{H^k}
	+\|u_0\|_{W^{3,1}}
	+\|\theta_0\|_{W^{3,1}}
	\le\varepsilon,
	\end{equation}
	the system \eqref{equation1} admits a unique global-in-time solution
	$(u,\theta)$ satisfying 
	\begin{equation*}
		\Vert u(t) \Vert_{H^k }^2 + \Vert \theta(t) \Vert_{H^k}^2 + \int_{0}^{t} (\Vert \partial_1 u(\tau) \Vert_{H^k}^2 + \Vert \partial_1 \theta(\tau) \Vert_{H^k}^2) d \tau  \leq C_0^2\varepsilon^2,
	\end{equation*}
for every $t\ge0$. Moreover, the solution satisfies the following decay estimates:
	\begin{align*}
		\Vert (u_1| \theta)(t) \Vert_{L^{2}} &\leq C_0 \varepsilon (1+t)^{-\frac{1}{8} - \eta }, \\ 
		\Vert u_2(t) \Vert_{L^{2}} &\leq C_0 \varepsilon (1+t)^{-\frac{1}{4} }, \\
		\Vert  (\partial_1 u_1| \partial_1 \theta)(t) \Vert_{L^{2}} &\leq C_0 \varepsilon (1+t)^{-\frac{5}{8} - \eta}, \\ 
		\Vert  (\partial_2 u_1| \partial_2 \theta)(t) \Vert_{L^{2}} &\leq C_0 \varepsilon (1+t)^{-\frac{1}{8} - \eta}, \\ 
		\Vert \partial_1 u_2(t) \Vert_{L^{2}} &\leq C_0 \varepsilon (1+t)^{-\frac{3}{4}}, \\
		\Vert  (u_1| \theta)(t) \Vert_{L^{4}} &\leq C_0 \varepsilon (1+t)^{-\frac{1}{2}}, \\
		\Vert u_2(t) \Vert_{L^{4}} &\leq C_0 \varepsilon (1+t)^{-\frac{7}{8}- \eta + \delta}, \\
		\Vert  (\partial_1u| \partial_1\theta)(t) \Vert_{L^{4}} &\leq C_0 \varepsilon (1+t)^{-1 +  \delta},
	\end{align*}
	where the parameters $\eta = \frac{1}{120}$,  $\delta = 10^{-5}$ and $C_0$ is a universal positive constant.
\end{thm}

Theorem~\ref{main2} establishes the global nonlinear stability of the hydrostatic equilibrium despite the presence of only horizontal dissipation. It shows that the hidden wave structure generated by the buoyancy coupling fundamentally changes the long-time dynamics of the system. Although the available dissipation acts only in the horizontal direction and is by itself insufficient to prevent the rapid growth predicted by standard energy estimates, the dispersive decay associated with the resulting internal gravity waves supplies the missing time integrability needed to close the nonlinear energy estimates. Consequently, the weak anisotropic dissipation, when combined with this dispersive stabilization, is sufficient to recover global stability. Moreover, the theorem provides explicit componentwise decay rates for the velocity and temperature that reflect the anisotropic nature of the underlying stabilization mechanism.

\vskip .1in 
The proof is based on the philosophy that the perturbation system should be viewed as a coupled dispersive-dissipative system rather than as the Navier--Stokes equations with an external forcing. The essential difficulty is that the horizontal dissipation alone is far too weak to control the nonlinear evolution. Instead, the missing decay is recovered from the hidden wave structure generated by the coupling between the velocity and temperature equations.

\vskip .1in 
The first step of the proof is to uncover the hidden dispersive structure of the perturbation system. Rather than analyzing the velocity and temperature equations separately, we combine them into the vector-valued unknown
\[
v=(u_1,u_2,\theta)^T
\]
and rewrite \eqref{equation1} as an evolution equation for $v$. After applying an extended Helmholtz projection, the linearized system takes the form
\[
\partial_t v+Lv+\widetilde{\mathbb P}J\widetilde{\mathbb P}v=0,
\]
where $L=\partial_1^2I$ represents the horizontal dissipation and the skew-adjoint operator $\widetilde{\mathbb P}J\widetilde{\mathbb P}$ encodes the coupling between the velocity and temperature.

\vskip .1in
A direct Fourier analysis shows that the operator $-\widetilde{\mathbb P}J\widetilde{\mathbb P}$ possesses the eigenvalues
	\begin{equation*}
	\Bigg\{i \frac{| \xi_1 |}{| \xi |}\ , -i \frac{| \xi_1 | }{| \xi |  }\ ,0 \Bigg\}.
\end{equation*}
Consequently, the linear solution admits the representation
	\begin{align*}
	T(t) v_0(x) = &  e^{t\partial_1^2}e^{t \mathcal{R}_1 }\mathcal{F}^{-1} \Big[\langle \hat{v}_0,\vec{a}_1(\xi)   \rangle_{\mathbb{C}^3} \vec{a}_1(\xi)\Big] (x) \\
	&+ e^{t\partial_1^2}e^{-t \mathcal{R}_1} \mathcal{F}^{-1} \Big[\langle \hat{v}_0,\vec{a}_2(\xi)   \rangle_{\mathbb{C}^3} \vec{a}_2(\xi)\Big] (x)\\
	&+ e^{t\partial_1^2} \mathcal{F}^{-1} \Big[\langle \hat{v}_0,\vec{a}_3(\xi)   \rangle_{\mathbb{C}^3} \vec{a}_3(\xi)\Big] (x),
\end{align*}
where $\mathcal R_1$ denotes the Fourier multiplier with symbol
$i|\xi_1|/|\xi|$.
The factors $e^{\pm t\mathcal R_1}$ describe the propagation of internal gravity waves, while the horizontal heat semigroup $e^{t\partial_1^2}$ provides the available anisotropic dissipation. Their interaction forms the fundamental dispersive--dissipative mechanism exploited throughout the paper.

\vskip .1in
Using the divergence-free condition, the non-propagating mode corresponding
to the zero eigenvalue vanishes identically. Consequently, both the linear
and nonlinear solutions can be represented entirely in terms of the two
oscillatory wave modes. More precisely, Duhamel's principle yields the
representation
\begin{equation}\label{Intro-Duhamel}
	v(t)
	=
	\sum_{\sigma\in\{\pm\}}
	S_\sigma(t)v_0
	-
	\int_0^t
	S_\sigma(t-s)\mathcal{N}(v(s))\,ds,
\end{equation}
where
\[
S_\sigma(t)
=
e^{t\partial_1^2}
e^{\sigma t\mathcal{R}_1}
Q_\sigma,
\]
$Q_\sigma$ denotes the projection onto the eigenspace corresponding to the
eigenvalue
$\sigma i\frac{|\xi_1|}{|\xi|}$,
and $\mathcal{N}(v)$ denotes the nonlinear terms.
The representation \eqref{Intro-Duhamel} reveals that the entire nonlinear
dynamics is governed by the oscillatory semigroups
$\{S_\pm(t)\}_{t\ge0}$.
The heart of the analysis is therefore to establish sharp dispersive decay
estimates for these semigroups and to combine them with the available
horizontal dissipation in the nonlinear energy estimates. The resulting
dispersive decay provides precisely the additional stabilization needed to
compensate for the lack of full dissipation and ultimately yields global
nonlinear stability.

\vskip .1in 
The proof then proceeds in two main steps. First, in Section~\ref{section3}, we establish the linear dispersive estimates for the oscillatory-dissipative semigroups $S_\pm(t)$. The key point is that the phase
\[
\pm\frac{|\xi_1|}{|\xi|}
\]
is anisotropic and degenerate in certain frequency regions. We therefore decompose the frequency space according to the relative sizes of $|\xi_1|$ and $|\xi_2|$. In the region where the oscillatory phase is non-degenerate, stationary phase yields dispersive decay. In the complementary region, where the phase degenerates, the horizontal heat kernel $e^{t\partial_1^2}$ provides the necessary decay. This combination of oscillation and horizontal dissipation leads to the crucial $L^4$ decay estimate for the linear flow.

\vskip .1in 
Second, in Section~\ref{section4}, we combine these linear decay estimates with refined nonlinear energy estimates through a bootstrap argument. The bootstrap assumptions include both uniform high-order Sobolev bounds and anisotropic decay estimates for the components of $(u,\theta)$ and their horizontal and vertical derivatives. Under these assumptions, the dispersive estimates obtained in Section~3 allow us to control the Duhamel terms and improve all bootstrap bounds by a factor of one half. The standard continuity argument then extends the solution globally in time and yields the decay estimates stated in Theorem~\ref{main2}.

\vskip .1in 
To place Theorem \ref{main2} into perspective, we briefly review several closely 
related developments on the stability of Boussinesq systems with partial dissipation. 
A first important observation is that the stabilizing mechanism exploited in this 
paper is genuinely tied to the stratified background: near the trivial equilibrium 
$(U,\Theta)=(0,0)$, the buoyancy force does not stabilize the flow and may in fact 
drive the growth of solutions. Brandolese and Schonbek \cite{BS} showed that, for 
the viscous Boussinesq system without thermal diffusion, the energy of the velocity 
field may grow in time even for small and well-localized data, in sharp contrast 
with the decay theory for the Navier--Stokes equations. Kukavica and Wang \cite{KW} 
investigated the long-time behavior of the two-dimensional Boussinesq equations 
with zero diffusivity and obtained quantitative upper bounds on the growth of 
Sobolev norms of the solution. More recently, Kiselev, Park and Yao \cite{KPY} 
proved the formation of small scales for the two-dimensional Boussinesq equations, 
establishing infinite-in-time growth of higher-order norms of the temperature. 
These results demonstrate that no decay or even uniform boundedness can be expected 
near the trivial steady state, and that the stability obtained in this paper is a 
manifestation of the wave structure generated by the hydrostatic stratification.

\vskip .1in 
For perturbations near hydrostatic equilibria, where the coupling between the 
velocity and the temperature provides an additional stabilizing mechanism, the 
picture is fundamentally different. Doering, Wu, Zhao and Zheng \cite{CJKX} 
initiated the rigorous study of the stability and large-time behavior of the 
hydrostatic equilibrium for the two-dimensional Boussinesq equations with only 
kinematic dissipation (and no thermal diffusion), working on a bounded domain 
with stress-free boundary conditions. In a subsequent paper, Tao, Wu, Zhao and 
Zheng \cite{TWZZ} resolved several questions left open in \cite{CJKX}, in 
particular establishing the velocity decay rates and identifying the eventual 
temperature profile, this time on the periodic domain. Stability has also been 
established when the viscous dissipation is replaced by weaker damping mechanisms. 
Wan \cite{Wan} proved the global well-posedness of the two-dimensional Boussinesq 
equations with a velocity damping term, Castro, C\'ordoba and Lear \cite{ADD} 
established the asymptotic stability of stratified solutions for the damped 
system, and Kim and Lee \cite{KimLee} further developed the stability theory 
for the stratified Boussinesq equations with velocity damping.

\vskip .1in 
The dispersive stabilization generated by the buoyancy coupling has subsequently 
been investigated in a variety of partially dissipative Boussinesq models. Ben 
Said, Pandey and Wu \cite{BPW} showed that the temperature coupling alone exerts 
a stabilizing effect on a class of buoyancy-driven flows. Lai, Wu and Zhong 
\cite{LJ} established global stability and large-time behavior for the 
two-dimensional Boussinesq equations with partial (degenerate) dissipation, 
exploiting the damped wave structure induced by the coupling, while Lai, Wu, Xu, 
Zhang and Zhong \cite{LWXZZ} derived optimal decay estimates. Adhikari, Ben Said, 
Pandey and Wu \cite{ABPW} proved global stability for the two-dimensional 
Boussinesq system with horizontal dissipation and vertical thermal diffusion, 
demonstrating how anisotropic dissipation and buoyancy-induced dispersion 
complement each other. Dong and Sun \cite{DCJ,DCJX} established the asymptotic 
stability of the two- and three-dimensional Boussinesq equations without thermal 
conduction. Subsequent works by Shang and Xu \cite{SX}, Wu and Zhang \cite{WZ}, 
Ji, Yan, and Wu \cite{JYW}, Kang, Lee, and Nguyen \cite{KLN}, Fujii and Li 
\cite{FujiiLi} and Lin, Wu and Zheng \cite{LinWuZheng} substantially advanced 
the theory by establishing stability and large-time behavior for anisotropic 
Boussinesq systems with partial dissipation near hydrostatic equilibrium. We 
refer the reader to these works and the references therein for further 
developments.

\vskip .1in 
It is worth emphasizing that the stability problem depends crucially on the underlying spatial domain. For the mixed domain $\Omega=\mathbb{T}\times\mathbb{R}$, Dong, Wu, Xu and Zhu \cite{DWXZ} established the global stability and exponential decay of the two-dimensional anisotropic Boussinesq equations with horizontal dissipation. Their analysis relies on decomposing the solution into its vertical average and oscillatory components and exploiting a Poincar\'e inequality for the oscillatory part. Such an approach is unavailable in the whole-space setting $\mathbb{R}^2$, where no distinguished zero Fourier mode exists and no Poincar\'e-type inequality is available. Consequently, the dispersive estimates developed in the present paper must completely replace the stabilizing mechanism provided by the periodic structure. Compared with the previous works, the present paper treats a substantially more delicate regime in which the available dissipation is purely horizontal, the velocity field is considerably more singular, and global stability is achieved almost entirely through the dispersive decay generated by the buoyancy coupling.

\textbf{Notations.}
In the $d$-dimensional Euclidean space $\mathbb{R}^d$, $L^p$ denotes the standard Lebesgue
space. The Sobolev space $W^{k,p}$ consists of functions whose weak derivatives up to order $k$ belong to $L^p$, while the homogeneous Sobolev space $\dot{W}^{k,p}$ is defined by requiring only the highest-order weak derivatives to belong to $L^p$. When $p=2$, $W^{k,p}$ and $\dot{W}^{k,p}$ are also denoted as ${H}^{k}$ and $\dot{H}^{k}$ correspondingly. Given two normed spaces $N_1$ and $N_2$ on measure spaces $(X_1,\mu_1)$ and
$(X_2,\mu_2)$, respectively, the anisotropic space
$N_2(X_2,\mu_2)\,N_1(X_1,\mu_1)$ is defined by first taking the $N_1$-norm in the $X_1$
variable and then the $N_2$-norm in the $X_2$ variable. More precisely, for a function
$f$ defined on the product space $X_1\times X_2$, we set
\[
\|f\|_{N_2(X_2,\mu_2)\,N_1(X_1,\mu_1)}
\triangleq \bigl\|\,\|f(x_1,x_2)\|_{N_1(X_1,\mu_1)}\,\bigr\|_{N_2(X_2,\mu_2)}.
\]
 Given two functions $f$ and $g$, we interpret the symbol $\Vert(f|g) \Vert_{N(X,\mu)}$ as follows:
\begin{equation*}
	\Vert(f|g) \Vert_{N(X,\mu)} \triangleq \Vert f \Vert_{N(X,\mu)} + \Vert g \Vert_{N(X,\mu)}.
\end{equation*}
The operator \(\Lambda_i^{s}\) denotes the fractional derivative in the
\(x_i\)-direction, defined via the Fourier transform by
\[
\widehat{\Lambda_i^{s} f}(\xi)
= |\xi_i|^{s}\widehat{f}(\xi),
\qquad 0\le s \le 1.
\]
The notation $a \lesssim b$ means that there exists  a uniform constant $C$ (which may be different in each occurrence), such that $a \leq Cb$. 

\vskip .1in
The rest of the paper is organized as follows. In Section~\ref{section2}, we symmetrize the linearized system and derive the explicit semigroup representation of the solution. In Section~\ref{section3}, we collect several useful lemmas and establish the key $L^4$ dispersive estimate for the linear part. In Section~\ref{section4}, we carry out the energy estimates and reduce the proof of Theorem~\ref{main2} to a family of decay estimates. Finally, in Section~\ref{section5}, we prove these decay estimates and complete the proof of Theorem~\ref{main2}.

\section{Symmetrization of Linear System}
\label{section2}

In this section, we perform the symmetrization of the linear system. First, we introduce a new unknown function which combines the velocity field with the thermal disturbance  
	\begin{equation*}
	v \triangleq  (u,\theta)^{T} = (u_1,u_2,\theta)^{T}.
	\end{equation*}
	Set 
	$$
	L \triangleq 
	\left(   \begin{array}{ccc}
	\partial_1^2 & 0 & 0 \\
	0 & \partial_1^2 & 0 \\
	0 & 0 & \partial_1^2
	\end{array} \right), \ \ 
	J \triangleq \left(   \begin{array}{ccc}
	0 & 0 & 0 \\
	0 & 0 & -1 \\
	0 & 1 & 0
	\end{array} \right),
	$$
	and $\widetilde{\nabla} \triangleq (\partial_1,\partial_2,0)^T$. Then the original system can be written as 
	\begin{align}\label{equation111}
	\left \{
	\begin{aligned} 
	& \partial_t v  + Lv  + Jv + (v \cdot \widetilde{\nabla}) v  + \widetilde{\nabla}p=    0,\\
	& \widetilde{\nabla} \cdot v =0    ,\\
	& v(0,x)= v_0(x).
	\end{aligned}
	\right.
	\end{align}
	Next, let $\widetilde{\mathbb{P}}$ be the extended Helmholtz projection, which is defined by
	$$
	\widetilde{\mathbb{P}} \triangleq
	\left(   \begin{array}{c|c}
	(\delta_{jk} + R_jR_k)_{1 \leq j,k \leq 2} & 0  \\
	\hline
	0 & 1 \\
	\end{array} \right).
	$$
Here $R_j$ is the Riesz transform defined by the Fourier multiplier $i\frac{\xi_j}{|\xi|}$. Applying the extended Helmholtz projection $\widetilde{\mathbb{P}}$ to (\ref{equation111}) gives the following evolution equation:
	\begin{align}\label{The system equation}
	\left \{
	\begin{aligned} 
	& \partial_t v  + Lv  + \widetilde{\mathbb{P}}J\widetilde{\mathbb{P}}v + \widetilde{\mathbb{P}}(v \cdot \widetilde{\nabla}) v  =    0,\\
	& \widetilde{\nabla} \cdot v =0    ,\\
	& v(0,x)= v_0(x).
	\end{aligned}
	\right.
	\end{align}
A direct computation shows that the eigenvalues of $-\tilde{\mathbb{P}} J \tilde{\mathbb{P}}$ are 
	\begin{equation*}
	\Bigg\{i \frac{| \xi_1 |}{| \xi |}\ , -i \frac{| \xi_1 | }{| \xi |  }\ ,0 \Bigg\},
	\end{equation*}
and	the corresponding eigenvectors $\{\vec{a}_j\}_{j=1}^3$ are given by 
	\begin{align*}
	&\vec{a}_1(\xi) \triangleq \frac{1}{\sqrt{2} |\xi_1 | |\xi |} (i\xi_1 \xi_2,-i\xi_1^2,|\xi_1 | |\xi |)^T, \\
	&\vec{a}_2(\xi) \triangleq \frac{1}{\sqrt{2} |\xi_1 | |\xi |} (-i\xi_1 \xi_2,i\xi_1^2,|\xi_1 | |\xi |)^T,\\
	&\vec{a}_3(\xi) \triangleq \frac{1}{|\xi | } (\xi_1 ,\xi_2,0)^T.
	\end{align*}
Straightforward computation shows that for $\xi_1 \neq 0$, the  eigenvectors $\{\vec{a}_j\}_{j=1}^3$ define an orthogonal basis in $\mathbb{C}^3$,
which is the 3-dimensional complex Euclidean space and the inner product is taken  as $\langle f , g \rangle_{\mathbb{C}^3} = \sum\limits_{j=1}^{3} f_j \bar{g}_j$.\\
So the semigroup $\{T (t)\}_{t \geq 0} $ generated by the linear operator $L+ \widetilde{\mathbb{P}} J \widetilde{\mathbb{P}}$ can be written explicitly by
	\begin{align*}
	T(t) v_0(x) = &  e^{t\partial_1^2}e^{t \mathcal{R}_1 }\mathcal{F}^{-1} \Big[\langle \hat{v}_0,\vec{a}_1(\xi)   \rangle_{\mathbb{C}^3} \vec{a}_1(\xi)\Big] (x) \\
	&+ e^{t\partial_1^2}e^{-t \mathcal{R}_1} \mathcal{F}^{-1} \Big[\langle \hat{v}_0,\vec{a}_2(\xi)   \rangle_{\mathbb{C}^3} \vec{a}_2(\xi)\Big] (x)\\
	&+ e^{t\partial_1^2} \mathcal{F}^{-1} \Big[\langle \hat{v}_0,\vec{a}_3(\xi)   \rangle_{\mathbb{C}^3} \vec{a}_3(\xi)\Big] (x),
	\end{align*}
where $\mathcal{R}_j$ denote the Fourier multiplier with symbol $i\frac{|\xi_j|}{|\xi|}.$
		By the divergence-free condition of $u_0$ and  the extended Helmholtz projection of the velocity $u \cdot \nabla u$ onto the divergence-free vector field, we  get that 
	\begin{align*}
	\langle \hat{v}_0,\vec{a}_3(\xi)   \rangle_{\mathbb{C}^3}  = 0, \\
	\langle \mathcal{F}{[\widetilde{\mathbb{P}}(v \cdot \tilde{\nabla}) v ] },\vec{a}_3(\xi)  \rangle_{\mathbb{C}^3} =0.
	\end{align*}
	Thus by  Duhamel principle, the solution of the nonlinear equation \eqref{The system equation} can be written as:
	\begin{align*}
	v(t,x) = &  e^{t\partial_1^2}e^{t\mathcal{R}_1 } \mathcal{F}^{-1} \Big[\langle \hat{v_0}, \vec{a}_1(\xi)   \rangle_{\mathbb{C}^3}   \vec{a}_1(\xi) \Big] (x) \\
	&+ e^{t\partial_1^2}e^{-t\mathcal{R}_1 } \mathcal{F}^{-1} \Big[\langle \hat{v_0},\vec{a}_2(\xi)   \rangle_{\mathbb{C}^3}   \vec{a}_2(\xi)\Big] (x) \\
	& + \int_{0}^{t } e^{(t-s)\partial_1^2}e^{(t-s)\mathcal{R}_1 } \mathcal{F}^{-1} \Big[\langle \mathcal{F}{[\widetilde{\mathbb{P}}(v \cdot \widetilde{\nabla}) v ] },\vec{a}_1(\xi)   \rangle_{\mathbb{C}^3}   \vec{a}_1(\xi)\Big] (x) ds \\
	& + \int_{0}^{t } e^{(t-s)\partial_1^2}e^{(t-s)\mathcal{R}_1 } \mathcal{F}^{-1} \Big[\langle \mathcal{F}{[\widetilde{\mathbb{P}}(v \cdot \widetilde{\nabla}) v ] },\vec{a}_2(\xi)   \rangle_{\mathbb{C}^3}   \vec{a}_2(\xi)\Big] (x) ds.
	\end{align*}
	Using  the explicit expressions of the eigenvectors, we can  write the solution as:
	\begin{equation}\label{formula of solution1}
		\begin{aligned}
	u_1 &= \sum_{\sigma \in \{ \pm \}}e^{t\partial_1^2}e^{\sigma t \mathcal{R}_1 } Q^{u_1}_{\sigma} (u_0,\theta_0) + \int_{0}^{t}e^{(t-s)\partial_1^2}e^{\sigma (t-s) \mathcal{R}_1 } Q^{u_1}_{\sigma} (\mathbb{P}(u \cdot \nabla u), u\cdot \nabla \theta ) \ ds , 
	\end{aligned}
	\end{equation}
		\begin{equation}\label{formula of solution2}
		\begin{aligned}
			u_2 &= \sum_{\sigma \in \{ \pm \}}e^{t\partial_1^2}e^{\sigma t \mathcal{R}_1 } Q^{u_2}_{\sigma} (u_0,\theta_0) + \int_{0}^{t}e^{(t-s)\partial_1^2}e^{\sigma (t-s) \mathcal{R}_1 } Q^{u_2}_{\sigma} (\mathbb{P}(u \cdot \nabla u), u\cdot \nabla \theta  ) \ ds,
		\end{aligned}
	\end{equation}
		\begin{equation}\label{formula of solution3}
		\begin{aligned}
			\theta &= \sum_{\sigma \in \{ \pm \}}e^{t\partial_1^2}e^{\sigma t \mathcal{R}_1 } Q^{\theta}_{\sigma} (u_0,\theta_0) + \int_{0}^{t}e^{(t-s)\partial_1^2}e^{\sigma (t-s) \mathcal{R}_1 } Q^{\theta}_{\sigma} (\mathbb{P}(u \cdot \nabla u), u\cdot \nabla \theta ) \ ds,
		\end{aligned}
	\end{equation}
	where 
	\begin{equation}\label{formula of solution4}
	\begin{aligned}
	&Q^{u_1}_{\pm} (u,\theta) \triangleq \frac{1}{2} \Big\{ - \frac{\partial_2^2}{\Lambda ^2} u_1 + \frac{\partial_1}{\Lambda }\frac{\partial_2}{\Lambda }u_2  \mp i \frac{\partial_1}{\Lambda_1 } \frac{\partial_2}{\Lambda } \theta \Big\}, 
	\end{aligned}
	\end{equation}
		\begin{equation}\label{formula of solution5}
		\begin{aligned}
			&Q^{u_2}_{\pm} (u,\theta) \triangleq \frac{1}{2} \Big\{  \frac{\partial_2}{\Lambda } \frac{\partial_1}{\Lambda } u_1 + \frac{\Lambda_1^2}{\Lambda ^2}u_2  \pm i \frac{\partial_1}{\Lambda_1 } \frac{\partial_1}{\Lambda } \theta \Big\}, 
		\end{aligned}
	\end{equation}
		\begin{equation}\label{formula of solution6}
		\begin{aligned}
			&Q^{\theta}_{\pm} (u,\theta) \triangleq \pm \frac{i}{2} \Big\{  \frac{\partial_1}{\Lambda_1} \frac{\partial_2}{\Lambda } u_1 - \frac{\partial_1}{\Lambda_1}\frac{\partial_1}{\Lambda }u_2  \Big\} + \frac{1}{2} \theta.
		\end{aligned}
	\end{equation}

\section{Tool Lemmas and Kernel Function Estimate}
\label{section3}

This section makes several preparations. First of all, we introduce the Littlewood-Paley decomposition (\ref{LP decomposition}). In the next subsection, we will give several useful inequalities which will be frequently used later. After introducing these useful inequalities, we shall present  the $L^4$ estimate for the linear system.

	\subsection{Littlewood-Paley decomposition}
	 Choose a radial function $\varphi \in \mathcal{S}(\mathbb{R}^2)$ supported in $ \{\xi \in \mathbb{R}^2, \frac{3}{4} \leq |\xi| \leq \frac{8}{3}\}$ such that
	\begin{equation*}
	\varphi(\xi) = 1, \ \ \mathrm{for}\ \  1 \leq |\xi | \leq \frac{4}{3},
	\end{equation*}
	and 
	\begin{equation*}
	\sum_{j \in \mathbb{Z}} \varphi(2^{-j}\xi ) = 1   \ \ \mathrm{for \  all } \ \ \xi \neq 0.
	\end{equation*}
For $f \in \mathcal{S}'_h(\mathbb{R}^2)$, the frequency localization operator $\Delta_j$ is defined by 
\begin{equation*}
\Delta_j f = \varphi(2^{-j}D)f.
\end{equation*}
	By localizing $f$ in frequency, we can decompose $f$ as
	\begin{equation}\label{LP decomposition}
	f = \sum_{j \in \mathbb{Z}}  \Delta_j f .
	\end{equation}
	\subsection{Some Tool Lemmas}
	 The following Bernstein’s lemma which is related to the Littlewood-paley decomposition will be repeatedly used throughout this paper.
	\begin{lem}\label{Bernstein}
		Let $\mathcal{C}$ be an annulus and $\mathcal{B}$ a ball. A constant $C$ exists such that for any nonnegative integer $k$, any couple $(p,q)$ in $[1,\infty]^2$ with $q \geq p \geq 1$, and any function $u$ in $L^p$, we have 
		\begin{align*}
		\mathrm{Supp} \ \hat{u} \subset \lambda \mathcal{B}  &\Longrightarrow \Vert D^k u \Vert_{ L^{q} } \triangleq  \sup\limits_{| \alpha | = k} \Vert \partial^\alpha u \Vert_{ L^{q} } \leq C^{k+1} \lambda^{k+d(\frac{1}{p} - \frac{1}{q})} \Vert u \Vert_{L^p} ,\\
		\mathrm{Supp} \ \hat{u} \subset \lambda \mathcal{C} &\Longrightarrow C^{-k-1} \lambda^k \Vert u \Vert_{ L^{p} } \leq \Vert D^k u \Vert_{ L^{p} } \leq C^{k+1} \lambda^k \Vert u \Vert_{ L^{p} }.
		\end{align*}
	\end{lem}
We state in the following lemma the fractional Leibniz rule and relevant commutator estimates. These inequalities are commonly named the Kato-Ponce inequalities, which  are of great utility for product estimates and will be met while dealing with nonlinear terms. The proofs of these inequalities are given in \cite{LiDong}. 
	\begin{lem}\label{Kato-Ponce inequality}
Let $f,g \in \mathcal{S}({\mathbb{R}}^d)$, let $s >0$, then there holds:
\begin{align*}
&\Vert D^s (fg) \Vert_{L^p} \leq C_{s,d,p_1,p_2,p_3,p_4}  (\Vert D^s f \Vert_{L^{p_1}} \Vert g \Vert_{L^{p_2}} + \Vert f \Vert_{L^{p_3}} \Vert D^s g \Vert_{L^{p_4}}), \\
& \Vert D^s (fg) - f D^s g \Vert_{L^p} \leq C_{s,d,p_1,p_2,p_3,p_4}  (\Vert D^s f \Vert_{L^{p_1}} \Vert g \Vert_{L^{p_2}} + \Vert D f \Vert_{L^{p_3}} \Vert D^{s-1}g \Vert_{L^{p_4}}), 
\end{align*}
where $$\frac{1}{p} = \frac{1}{p_1} + \frac{1}{p_2} = \frac{1}{p_3} + \frac{1}{p_4},\  1 < p < \infty, \ 1 < p_1, p_2, p_3, p_4 \leq \infty$$ and $C_{s,d,p_1,p_2,p_3,p_4}$ is a constant depending only on $s,d,p_1,p_2,p_3,p_4$.
\end{lem}
	The following lemma provides the decay estimate for a convolution-type integral which will be used when we investigate the decay for the nonlinear Duhamel integral. Its proof can be found in \cite{MYZ}.
	\begin{lem}\label{decay lem}
		Assume $\alpha \ge 1$, $\beta >1 $. Then, for some constant $C> 0$ and $t>1$, we have
		\begin{align}\label{decay first}
		&\int_0^{t-1}(t-\tau)^{-\alpha}(1+\tau)^{-\beta}\:d\tau\leq
		\begin{cases}
		C (1+t)^{-\alpha},\ &\alpha \le \beta,\\
		C (1+t)^{-\beta},\ &\alpha > \beta.\\
		\end{cases}
		\end{align}
		Assume $\alpha < 1$. Then, for some constant $C> 0$, we have
		\begin{align}\label{decay second}
		&\int_0^{t}(t-\tau)^{-\alpha}(1+\tau)^{-\beta}\:d\tau\leq
		\begin{cases}
		C(1+t)^{-\alpha},\ &\beta >1,\\
		C(1+t)^{-\alpha} \ln(1+t),\ &\beta=1,\\
		C(1+t)^{1-\alpha-\beta},\ &\beta<1.
		\end{cases}
		\end{align}
	\end{lem}
	The  heat semigroup estimate is frequently used in the following text,  the proof of which can be found in \cite{MYZ}.
	\begin{lem}\label{HEAT DECAY LEMMA}
		Let $\sigma \ge 0$, $\alpha>0$, $\nu>0$, $1 \leq p \leq q \leq \infty$. Then 
		\begin{align*}
		\|\Lambda^\sigma e^{-\nu(-\Delta)^\alpha t} f\|_{L^q\left(\mathbb{R}^d\right)} \le C t^{-\frac{\sigma}{2 \alpha}-\frac{d}{2 \alpha}\left(\frac{1}{p}-\frac{1}{q}\right)}\|f\|_{L^p\left(\mathbb{R}^d\right)}.
		\end{align*}
	\end{lem}
The next lemma provides an inequality which reflects the action of the heat flow on spectrally supported functions. One may refer to \cite{C} for the detailed proof.
\begin{lem}\label{Localized heat kernel}
	There exist absolute positive constants c and C such that
	\begin{equation*}
	\Vert e^{ \Delta t} \Delta_j u \Vert_{L^p(\mathbb{R}^d)} \leq Ce^{-c2^{2j}t}	\Vert  \Delta_j u \Vert_{L^p(\mathbb{R}^d)} 
	\end{equation*}
	for all $j \in \mathbb{Z}$, $t > 0$. Moreover, for any $s >0$, there exists a positive constant $C = C(s)$ such that 
	\begin{equation*}
	\sum_{j \in \mathbb{Z}} 2^{sj} \Vert e^{ \Delta t} \Delta_j u \Vert_{L^p(\mathbb{R}^d)} \leq C t^{-\frac{s}{2}} \sup_{j \in \mathbb{Z}} \Vert  \Delta_j u \Vert_{L^p(\mathbb{R}^d)} .
	\end{equation*}
\end{lem}

\subsection{The linear estimate}
This subsection aims to establish linear estimates, specifically an $L^4$ estimate for the linear part that couples the dispersive and dissipative effects. Recall the representation formula for solutions to the linear part:
	\begin{align*}
u_1^{lin} & \triangleq \sum_{\sigma \in \{ \pm \}}e^{t\partial_1^2}e^{\sigma t \mathcal{R}_1 } Q^{u_1}_{\sigma} (u_0,\theta_0),  \quad
u_2^{lin}  \triangleq \sum_{\sigma \in \{ \pm \}}e^{t\partial_1^2}e^{\sigma t \mathcal{R}_1 } Q^{u_2}_{\sigma} (u_0,\theta_0) ,\\
\theta^{lin} & \triangleq \sum_{\sigma \in \{ \pm \}}e^{t\partial_1^2}e^{\sigma t \mathcal{R}_1 } Q^{\theta}_{\sigma} (u_0,\theta_0),
\end{align*}
where $Q^{u_1}$, $Q^{u_2}$ and $Q^{\theta}$ are defined in \eqref{formula of solution4}--\eqref{formula of solution6}.

By Young's convolution inequality, it is direct to obtain that the Riesz transform (when $d=1$, it is the Hilbert transform) is uniformly bounded in the sense that there exists a constant $C > 0$ independent of $j$ such that:
\begin{align}\label{Riesz operator}
\Vert R_j \Delta_j f\Vert_{L^{p}(\mathbb{R}^d)} \leq C \Vert \Delta_j f \Vert_{L^p(\mathbb{R}^d)} \ \ \mathrm{for \  all } \ 1 \leq p \leq \infty.
\end{align}
So the linear estimate is reduced to the decay estimate for the kernel $e^{(\partial_1^2 + \mathcal{R}_1)t}$. For simplicity, let us study the following linear equation, whose solution is generated by the semigroup $e^{(\partial_1^2 + \mathcal{R}_1)t}$,
	\begin{align}\label{linear equation}
\left \{
\begin{aligned} 
&  \partial_t f - \partial_1^2 f - \mathcal{R}_1 f =0, \ t > 0,\\
&f(x,0) = f_0(x).
\end{aligned}
\right.
\end{align}
Roughly speaking, we divide the frequency space into two parts, $ |\xi_1 | \leq | \xi_2| $ and $ |\xi_2 | \leq | \xi_1|$. In the first case, the kernel $e^{i\frac{\xi_1}{| \xi | } t }$ is non-degenerate, so it can make a better dispersive decay estimate; in the second case, the kernel $e^{i\frac{\xi_1}{| \xi | } t }$ is degenerate, it can only make a degenerate dispersive estimate, but the 1-D kernel $e^{-\xi_1^2t}$ behaves as the 2-D kernel $e^{-(\xi_1^2+ \xi_2^2)t}$.

In order to investigate the dispersive estimate, we use the stationary phase estimate of oscillatory integrals under small perturbations. Its proof can be found in \cite{Stein}.
\begin{lem}\label{stationary phase lemma}
	Let $d \geq 1$ be an integer. Let $\psi \in C_c^{\infty} (\mathbb{R}^d) $ and let $U$ be a neighborhood of supp $\psi$. Let $\Phi \in C^{\infty}(U;\mathbb{R}) $  satisfy
	\begin{equation*}
	rank(\nabla^2 \Phi (\xi) ) \geq k, \           \      \xi \in supp\ \{\psi\},
	\end{equation*}
	for some $k \in \{1, \dots d\}$. Then there exist two positive constants $\varepsilon = \varepsilon(d,\psi,\Phi)$ and $C = C(d,\psi ,\Phi)$ such that 
	\begin{equation*}
	\left|  \int_{R^d}    e^{i x \cdot \xi}e^{it (\Phi(\xi) + \Psi(\xi))}       \psi(\xi)  d\xi  \right| \leq C(1+\vert t \vert)^{-\frac{k}{2}},
	\end{equation*}
	for all $t \in \mathbb{R}, x \in \mathbb{R}^d$ and $\Psi \in C^{\infty}(U;\mathbb{R})$ satisfying $\Vert \Psi \Vert_{C^{d+3}(supp \{\psi\}) } \leq \varepsilon$.
\end{lem}
In the next lemma, we establish the localized $L^2$ estimate and $L^\infty$ estimate for the semigroup $e^{(\partial_1^2 + \mathcal{R}_1)t}$.
\begin{lem}\label{localized estimate}
For $f_0 \in L^1(\mathbb{R}^2) \bigcap L^2(\mathbb{R}^2)$,  we have
\begin{align}
& \  \Vert   e^{(\partial_1^2 + \mathcal{R}_1)t} \Delta_j f_0 \Vert_{L^{2}} \lesssim    \Vert  \Delta_j  f_0 \Vert_{L^{2}} \label{localized L2 estimate} ,
\\
& \ \Vert  e^{(\partial_1^2 + \mathcal{R}_1)t} \Delta_j f_0 \Vert_{L^{\infty}} 
\lesssim   (2^{2j} (1+t)^{-1} + 2^{2j} e^{-c2^{2j}t} (1+t)^{-\frac{1}{2}})\Vert  \Delta_j f_0 \Vert_{L^{1}} \label{localized Linfty estimate}.
\end{align}	
\end{lem}
\noindent Proof: 	According to lemma \ref{Bernstein} and Plancherel's Theorem, it is straightforward to obtain that 
	\begin{align*}
	 \Vert   e^{(\partial_1^2 + \mathcal{R}_1)t} \Delta_j f_0 \Vert_{L^{2}}   
	= &  \Vert  e^{ i \frac{| \xi_1 |}{|\xi|} t } e^{-\xi_1^2 t}\varphi(2^{-j}\xi)   \hat{f}_0 \Vert_{L^{2}} 
	\lesssim    \Vert  \Delta_j  f_0 \Vert_{L^{2}} .
	\end{align*}
Then we show the proof of (\ref{localized Linfty estimate}) by virtue of Lemma \ref{HEAT DECAY LEMMA}. For the kernel $e^{ i \frac{| \xi_1 |}{|\xi|} t }$, its phase function $\frac{| \xi_1 |}{|\xi|}$ is not smooth at $\xi_1 = 0$, so in order to use Lemma \ref{stationary phase lemma}, we write it as 
\begin{equation*}
e^{ i \frac{| \xi_1 |}{|\xi|} t } = \frac{1+sgn(\xi_1)}{2}e^{ i \frac{ \xi_1 }{|\xi|} t } +  \frac{1-sgn(\xi_1)}{2} e^{ -i \frac{ \xi_1 }{|\xi|} t },
\end{equation*}
where $sgn(\xi_1)$ is the sign function of $\xi_1$, which is the symbol of the Hilbert transform.
So by the boundedness of the Hilbert transform \eqref{Riesz operator},  we get 
\begin{equation*}
\Vert  e^{(\partial_1^2 + \mathcal{R}_1)t} \Delta_j f_0 \Vert_{L^{\infty}} \lesssim \Vert  e^{(\partial_1^2 + R_1)t} \Delta_j f_0 \Vert_{L^{\infty}}. 
\end{equation*}
A direct computation  gives the Hessian (in $\xi$) of $\phi \triangleq \frac{\xi_1}{| \xi |} $ as 
	\begin{equation*}
	H_{\phi} (\xi) = \vert \xi \vert^{-5} 
	\begin{pmatrix}
	-3 \xi_1 \xi_2^2 & \xi_2(2\xi_1^2 - \xi_2^2) \\
	\xi_2(2\xi_1^2 - \xi_2^2) & -\xi_1(\xi_1^2 - 2 \xi_2^2)
	\end{pmatrix}
	\end{equation*}
	with determinant $- \frac{\xi_2^2}{\vert \xi \vert^6}$.\\
 Since  the  kernel function $\phi$ has a degenerate point $\xi_2 = 0$, in order to achieve better  decay rate, we introduce a smooth bump function $\psi_1 : \mathbb{R} \to \mathbb{R}$  satisfying:
\begin{equation*}
\psi_1(\xi_1) = 1,\  \ \frac{1}{3} \leq | \xi_1 | \leq 3,  
\end{equation*} 
and 
\begin{equation*}
	\psi_1(\xi_1 ) = 0, \ \ |\xi_1 | \leq \frac{1}{4} \  \mathrm{or} \  | \xi_1 | \geq 4.
\end{equation*}
 Let $\psi_2(\xi_1) = 1- \psi_1(\xi_1)$. Then  it is easy to verify  that  
 \begin{align*}
 \mathrm{Supp}  \{\psi_1(\xi) \varphi(\xi)\} \subset \Big\{ \xi \in \mathbb{R}^2  \Big| \ \frac{3}{4} \leq |\xi| \leq \frac{8}{3},\ \frac{1}{4} \leq |\xi_1| \leq \frac{8}{3} \Big\}
 \end{align*}
 and 
  \begin{align*}
 \mathrm{Supp}  \{\psi_2(\xi) \varphi(\xi)\} \subset \Big\{ \xi \in \mathbb{R}^2  \Big| \ \frac{3}{4} \leq |\xi| \leq \frac{8}{3},\ |\xi_1| \leq \frac{1}{3},\ \frac{\sqrt{65}}{12} \leq |\xi_2| \leq \frac{8}{3} \Big\}.
 \end{align*}
 For the $L^\infty$ estimate of $\Delta_j f_0$, we use the Young's  inequality for convolution to obtain that
	\begin{align*}
	 \Vert  e^{(\partial_1^2 + {R}_1)t} \Delta_j f_0 \Vert_{L^{\infty}} 
	= & \Vert  e^{(\partial_1^2 + {R}_1)t} (\psi_1(2^{-j} D_1) + \psi_2(2^{-j} D_1) )\Delta_j f_0 \Vert_{L^{\infty}} \\
	\lesssim  &\Vert ( e^{ i \frac{\xi_1}{|\xi|} t }  \tilde{\varphi}(2^{-j}\xi) )^{\vee} \ast e^{\partial_1^2 t} \psi_1(2^{-j} D_1) \Delta_j  f_0 \Vert_{L^{\infty}} \\
	& + \Vert ( e^{ i \frac{\xi_1}{|\xi|} t }  {\psi}_2(2^{-j} \xi_1)  \tilde{\varphi}(2^{-j}\xi) )^{\vee} \ast e^{\partial_1^2 t}    \Delta_j f_0 \Vert_{L^{\infty}} \\
	\lesssim &  \Vert ( e^{ i \frac{\xi_1}{|\xi|} t }   \tilde{\varphi} (2^{-j}\xi) )^{\vee} \Vert_{L^{\infty}} \Vert e^{\partial_1^2 t}  \psi_1(2^{-j} D_1) \Delta_j f_0 \Vert_{L^{1}} \\
	& + \Vert ( e^{ i \frac{\xi_1}{|\xi|} t }  \psi_2(2^{-j} \xi)  \tilde{\varphi}(2^{-j}\xi) )^{\vee} \Vert_{L^{\infty}} \Vert   \Delta_j f_0 \Vert_{L^{1}}\\
	\triangleq & A_1 + A_2,
	\end{align*}
	where $\tilde{\varphi} $ is a  slightly “fattened” bump function (which equals 1 on the support of $\varphi $) satisfying 
	\begin{equation*}
	\tilde{\varphi} (\xi) = 0, \quad \quad | \xi | \leq \frac{1}{2} \quad \mathrm{or} \quad |\xi| \geq \frac{25}{9}.
	\end{equation*}
By the definition of $\tilde{\varphi}(\xi)$, we can obtain that
	  \begin{align*}
	\mathrm{Supp}  \{\psi_2(\xi) \tilde{\varphi}(\xi)\} \subset \Big\{ \xi \in \mathbb{R}^2  \Big| \ \frac{1}{2} \leq |\xi| \leq \frac{25}{9},\ |\xi_1| \leq \frac{1}{3},\ \frac{\sqrt{5}}{6} \leq |\xi_2| \leq \frac{25}{9} \Big\}.
	\end{align*}
		For the first term $A_1$, by scaling, we can obtain
	\begin{align*}
  \mathcal{F}^{-1} \Big[( e^{ i \frac{\xi_1}{|\xi|} t }   \tilde{\varphi}(2^{-j}\xi) ) \Big](x)  	&=   \int_{\mathbb{R}^2} e^{i x \cdot \xi } e^{i\frac{\xi_1}{  \vert \xi \vert }t} \tilde{\varphi}(2^{-j} \xi)  d\xi \\
	=  & 2^{2j} \int_{\vert \xi \vert \sim 1} e^{i 2^j x \cdot \xi } e^{i\frac{\xi_1}{  \vert \xi \vert }t}  \varphi( \xi)  d\xi \\ 
	= & 2^{2j}\int_{\vert \xi \vert \sim 1} e^{i 2^j x \cdot \xi } e^{ i\phi(\xi) t }   \varphi( \xi)  d\xi,
	\end{align*}
	by Lemma \ref{stationary phase lemma}, we have  
	\begin{align*}
	\Big| \int_{\vert \xi \vert \sim 1}  e^{i (2^j x) \cdot \xi }e^{ i\frac{\xi_1}{  \vert \xi \vert } t}   \varphi( \xi)  d\xi \Big| \lesssim    (1+t)^{-\frac{1}{2}}.
	\end{align*} 
The kernel function $\phi$ is degenerate, but we can  use the heat kernel to obtain the decay.	By the definition of cut-off operator $\psi_1$ and Lemma \ref{Localized heat kernel},  we have 
	\begin{align*}
	\Vert e^{\partial_1^2 t}  \psi_1(2^{-j} D_1) \Delta_j f_0 \Vert_{L^{1}} \lesssim e^{-c2^{2j}t} \Vert    \Delta_j f_0 \Vert_{L^{1}}.
	\end{align*}
As a result, by Lemma \ref{Localized heat kernel},
\begin{equation}\label{estimate for A1}
	\begin{aligned}
	A_1 & \lesssim 2^{2j} e^{-c2^{2j}t}(1+t)^{-\frac{1}{2}} .\\
	\end{aligned}
\end{equation}
For the second term,  by scaling 
	\begin{align*} 
	&(( e^{ i \frac{\xi_1}{|\xi|} t } \psi_2(2^{-j} \xi_1)  \tilde{\varphi}(2^{-j}\xi) )^{\vee}(x) \\
	= & \int_{\mathbb{R}^2} e^{i x \cdot \xi } e^{i\frac{\xi_1}{\vert \xi \vert }t} \psi_2(2^{-j} \xi_1) \tilde{\varphi}(2^{-j} \xi)  d\xi \\
	= & 2^{2j} \int_{\vert \xi \vert \sim 1,\  \vert \xi_2 \vert \sim 1} e^{i (2^j x) \cdot \xi } e^{ i \frac{\xi_1}{\vert \xi \vert } t }   \psi_2( \xi_1) \tilde{\varphi} ( \xi)  d\xi \\ 
	= & 2^{2j} \int_{\vert \xi \vert \sim 1,\  \vert \xi_2 \vert \sim 1}  e^{i (2^j x) \cdot \xi } e^{ i\phi(\xi) t }   \psi_2( \xi_1)\tilde{\varphi}( \xi)  d\xi.
	\end{align*}
Owing to the presence of $\psi_2$, the  kernel function $\phi$ is non-degenerate, so by Lemma \ref{stationary phase lemma}, we get
	\begin{align}\label{estimate for A2}
	\Big| \int_{\vert \xi \vert \sim 1,\  \vert \xi_2 \vert  \sim 1} e^{i 2^j x \cdot \xi } e^{ i\phi(\xi) t }   \varphi( \xi)  d\xi \Big| \lesssim (1+t)^{-1}.
	\end{align} 
This leads to
  \begin{align}
  	A_2 \lesssim 2^{2j} (1+t)^{-1} \Vert \Delta_j f_0 \Vert_{L^{1}}.
  \end{align}
Combining these two estimates \eqref{estimate for A1} and \eqref{estimate for A2} for $A_1$ and $A_2$ yields the desired bound \eqref{localized Linfty estimate}. 
	
	In order to obtain the $L^p \ (1 < p < \infty)$ estimate, we need to use the following Riesz-Th\"orin interpolation Theorem, the proof of which can be found in \cite{F}.
	\begin{lem}  \label{Riesz-Thorin}
		Suppose that (X, M, $\mu$) and (Y, N, $\nu$) are measure spaces and $p_0, p_1, q_0, q_1 \in [1,\infty]$. If $q_0 = q_1 = \infty$, suppose also that $\nu$ is semifinite. For $0 < t < 1$, define $p_t$ and $q_t$ by
		\begin{equation*}
		\frac{1}{p_t} = \frac{1-t}{p_0} + \frac{t}{p_1},\ \ \  \frac{1}{q_t} = \frac{1-t}{q_0} + \frac{t}{q_1}.
		\end{equation*}
		If  $T$ is a linear map from $L^{p_0}(\mu) + L^{p_1}(\mu)$ into $L^{q_0}(\nu) + L^{q_1}(\nu)$ such 
		that 
		\begin{equation*}
		\Vert Tf \Vert_{q_0} \leq M_0 \Vert f \Vert_{p_0} , \ \ \forall f \in L^{p_0}(\mu),
		\end{equation*}  
		and 
		\begin{equation*}
		\Vert Tf \Vert_{q_1} \leq M_1 \Vert f \Vert_{p_1}  , \ \ \forall f \in L^{p_1}(\mu),
		\end{equation*}
		then 
		\begin{equation*}
		\Vert Tf \Vert_{q_t} \leq M_0^{1-t}M_1^{t} \Vert f \Vert_{p_t} , \ \ \forall  f \in L^{p_t}(\mu).
		\end{equation*}
	\end{lem}
Using Lemma \ref{Riesz-Thorin} and the estimate in Lemma \ref{localized estimate}, we can obtain the following $L^4$ estimate for the linear equation \eqref{linear equation}.
	\begin{pro}\label{3} If $f_0 \in W^{1+\gamma,\frac{4}{3}}$, the solution $f$ of equation (\ref{linear equation}) has the following decay estimate:
		\begin{equation}\label{L4 decay estimate}
		\Vert  f \Vert_{L^{4}} \lesssim (1+t)^{-\frac{1}{2}} \Vert f_0 \Vert_{W^{1+\gamma,\frac{4}{3}}},
		\end{equation}
where the parameter $\gamma = 10^{-5}$.
	\end{pro}
	\noindent Proof : 
	By Lemma \ref{Riesz-Thorin} and Lemma \ref{localized estimate}, we get 
	\begin{equation}\label{localized L4 to L43}
	\Vert  e^{(\partial_1^2 + \mathcal{R}_1)t} \Delta_j f_0 \Vert_{L^{4}} \lesssim (2^j (1+t)^{-\frac{1}{2}} +2^j e^{-c2^{2j}t}(1+t)^{-\frac{1}{4}}) \Vert  \Delta_j f_0 \Vert_{L^{\frac{4}{3}}}.
	\end{equation}
	Therefore
	\begin{align*}
	\Vert f \Vert_{L^{4}} 
	\lesssim &\sum_{j \in \mathbb{Z}} \Vert  e^{(\partial_1^2 + \mathcal{R}_1)t} \Delta_j f_0 \Vert_{L^{4}} \\
	\lesssim & \sum_{j \in \mathbb{Z}} (2^{\frac{j}{2}} + 2^j) (1+t)^{-\frac{1}{2}} \Vert \Delta_j f_0 \Vert_{L^{\frac{4}{3}}}   \\
	\lesssim & (1+t)^{-\frac{1}{2}} (\Vert  f_0 \Vert_{L^{\frac{4}{3}}} +\Vert   \Lambda^{1+\gamma} f_0 \Vert_{L^{\frac{4}{3}}} ) \\
	\lesssim & (1+t)^{-\frac{1}{2}} \Vert f_0 \Vert_{W^{1+\gamma,\frac{4}{3}}}.
	\end{align*}
	Thus we complete the proof of (\ref{L4 decay estimate}).
	\section{Energy estimate and Proof of Theorem \ref{main2} }\label{section4}
In the next two sections, we prove Theorem \ref{main2}. Since the local well-posedness of (\ref{equation1}) can be obtained by a standard approach such as Friedrichs' method of cutoff in Fourier space, we focus on deriving the global bounds of $u$ and $\theta$. 

The argument is based on a bootstrap approach combined with refined energy estimates that yield global-in-time decay and stability.\\
\noindent\textbf{Bootstrap assumptions.}
We assume that there exists a constant $C_0 > 0$ such that for all $t \in [0,T )$, the solution satisfies the following uniform bounds: 
\begin{align}
\Vert (u| \theta)(t) \Vert_{H^k} &\leq C_0 \varepsilon,  \label{equation15}\\
\Vert \partial_1^\alpha(u_1| \theta)(t) \Vert_{L^{2}} &\leq C_0 \varepsilon (1+t)^{-\frac{1+4\alpha}{8} - \eta }, \\ 
\Vert  \partial_1^\alpha (u_1| \theta)(t) \Vert_{L^{4}} &\leq C_0 \varepsilon (1+t)^{-\frac{1+(1-\delta)\alpha}{2}}, \\
\Vert \partial_1^\alpha u_2(t) \Vert_{L^{2}} &\leq C_0 \varepsilon (1+t)^{-\frac{1+2\alpha}{4} }, \\
\Vert u_2(t) \Vert_{L^{4}} &\leq C_0 \varepsilon (1+t)^{-\frac{7}{8}- \eta + \delta}, \\
\Vert  \partial_1^\alpha \partial_2( u_1| \theta)(t) \Vert_{L^{2}} &\leq C_0 \varepsilon (1+t)^{-\frac{1+4\alpha}{8} - \eta}, \label{equation16}
\end{align}
where the parameters $\alpha \in \{0,1\}$, $\eta = \frac{1}{120}$ and $\delta = 10^{-5}$. \\
\noindent\textbf{Improved estimates.}
Based on the initial condition (\ref{Initial data condition}) and the above \emph{a priori} bounds,
we will establish the improved estimates
\begin{align}
\Vert (u| \theta)(t) \Vert_{H^k} &\leq \frac{1}{2}C_0 \varepsilon,   \label{IEu}\\
\Vert \partial_1^\alpha(u_1| \theta)(t) \Vert_{L^{2}} &\leq \frac{1}{2}C_0 \varepsilon (1+t)^{-\frac{1+4\alpha}{8} - \eta },\label{IDu1} \\ 
\Vert \partial_1^\alpha (u_1| \theta)(t) \Vert_{L^{4}} &\leq \frac{1}{2}C_0 \varepsilon (1+t)^{-\frac{1+(1-\delta)\alpha}{2}}, \label{IDuL4}\\
\Vert \partial_1^\alpha u_2(t) \Vert_{L^{2}} &\leq \frac{1}{2}C_0 \varepsilon (1+t)^{-\frac{1+2\alpha}{4} }, \label{IDu2}\\
\Vert u_2(t) \Vert_{L^{4}} &\leq \frac{1}{2}C_0 \varepsilon (1+t)^{-\frac{7}{8}- \eta + \delta}, \label{IDu2L4}\\
\Vert  \partial_1^\alpha \partial_2( u_1 | \theta)(t) \Vert_{L^{2}} &\leq \frac{1}{2}C_0 \varepsilon (1+t)^{-\frac{1+4\alpha}{8} - \eta}. \label{IDpartial2u}
\end{align}
Once (\ref{IEu})—(\ref{IDpartial2u}) are established, a standard bootstrap argument implies that the maximal existence time satisfies $T= \infty$, and hence all decay estimates hold globally in time.
\subsection{Estimate for $\Vert (u | \theta) \Vert_{H^k}$}
 In this section, we will verify (\ref{IEu}) first. By virtue of the norm equivalence
\begin{equation*}
\Vert (u| \theta) \Vert_{H^k} \sim \Vert (u | \theta) \Vert_{L^2} +  \Vert  (u | \theta) \Vert_{\dot{H}^k},
\end{equation*}
it suffices to bound the $L^2$ norm and $\dot{H}^{k}$ norm of $u$ and $ \theta$.
The $L^2$ energy estimate for $u$ and $\theta$ is standard,
\begin{equation*}
\Vert (u | \theta)(t) \Vert_{L^2}^2 + 2\int_{0}^{t} \Vert \partial_1 u(\tau) \Vert_{L^2}^2 d\tau + 2\int_{0}^{t}\Vert \partial_1 \theta(\tau) \Vert_{L^2}^2 d\tau = \Vert (u_0 | \theta_0) \Vert_{L^2}^2 .
\end{equation*}
Applying $\partial_m^k $ (with $m\in\{1,2\}$) to the equation of $u$ and  $\theta$,  then taking the $L^2$ inner product with $\partial_m^k u$  and $\partial_m^k \theta$ respectively. Adding the resulting estimates together and integrating from $0$ to $t$, we get
\begin{equation}\label{equation11}
\begin{aligned}
& \Vert (u| \theta)(t) \Vert_{H^k}^2 + 2\int_{0}^{t} \Vert \partial_1 u(\tau) \Vert_{H^k}^2 d \tau +2 \int_{0}^{t}\Vert \partial_1 \theta(\tau) \Vert_{H^k}^2 d\tau \\
= & \Vert (u_0| \theta_0) \Vert_{H^k}^2 -2\int_{0}^{t} (B_1 + B_2+ B_3 + B_4 )(\tau) d \tau,
\end{aligned}
\end{equation}
where 
\begin{align*}
& B_1 = (\partial_1^k (u \cdot \nabla u),\partial_1^k u )_{L^2} + (\partial_2^k (u \cdot \nabla u_2),\partial_2^k u_2 )_{L^2},\\
& B_2 = (\partial_1^k (u \cdot \nabla \theta),\partial_1^k \theta )_{L^2} , \\
& B_3 = (\partial_2^k (u \cdot \nabla \theta),\partial_2^k \theta )_{L^2},\\
& B_4 = (\partial_2^k (u \cdot \nabla u_1),\partial_2^k u_1 )_{L^2}.
\end{align*}
Now we estimate each term on the right-hand side of (\ref{equation11}). For 
\begin{align*}
B_1 = \sum_{i=1}^{2}\sum_{j = 1}^{k} C_{k}^{j} (\partial_1^j u_i \partial_i \partial_1^{k-j}  u, \partial_1^k u )_{L^2} + \sum_{i=1}^{2}\sum_{j = 1}^{k}  C_k^j (\partial_2^j u_i  \partial_i \partial_2^{k-j} u_2,\partial_2^k u_2 )_{L^2},
\end{align*}
by divergence-free condition of $u$,  H\"older's inequality, it is easy to derive that
\begin{equation*}
	| B_1 | \lesssim \Vert u \Vert_{H^k} \Vert \partial_1 u \Vert_{H^k}^2.
\end{equation*}
Similarly, $B_2$ can be estimated by analogy with $B_1$:
\begin{equation*}
	| B_2 | \lesssim \Vert ( u | \theta ) \Vert_{H^k} \Vert  (\partial_1 u | \partial_1 \theta)\Vert_{H^k}^2.
\end{equation*}
We now focus on the estimation of $B_3$. Since $u$ satisfies the divergence-free condition, we can get
\begin{align*}
\vert B_3  \vert  = & |(\partial_2^k (u \cdot \nabla \theta)-u \cdot \nabla \partial_2^k \theta ,\partial_2^k \theta )_{L^2}|  \\
\lesssim & |(\partial_2^k (u_1 \partial_1  \theta)- u_1\partial_1 \partial_2^k \theta,\partial_2^k \theta )_{L^2}| + |(\partial_2^k (u_2 \partial_2 \theta) - u_2 \partial_2^{k+1}\theta,\partial_2^k \theta)_{L^2}| \\
\triangleq & B_{31} + B_{32}.
\end{align*}
For $B_{31}$, by H\"older's inequality, Lemma \ref{Kato-Ponce inequality}  we can bound it by: 
\begin{align*}
B_{31} \lesssim   &\Vert (\partial_2^k (u_1 \partial_1  \theta)- u_1\partial_1 \partial_2^k \theta) \Vert_{L^2_{x_1}(L^2_{x_2})} \Vert \partial_2^k \theta \Vert_{L^2} \\
\lesssim & \Big( \Big\Vert \Vert \partial_2 u_1 \Vert_{L^4_{x_2}} \Vert \partial_1 \partial_2^{k-1}\theta \Vert_{L^4_{x_2}} \Big\Vert_{L^2_{x_1}}+  \Big\Vert \Vert \partial_2^k u_1 \Vert_{L^2_{x_2}} \Vert \partial_1 \theta \Vert_{L^{\infty}_{x_2}} \Big\Vert_{L^2_{x_1}} \Big)\Vert \partial_2^k \theta \Vert_{L^2}\\
\lesssim &  \Vert \partial_2 u_1 \Vert_{L^{\infty}_{x_1}(L^4_{x_2}) }\Vert \partial_1 \partial_2^{k-1}\theta \Vert_{L^2_{x_1}(L^4_{x_2})} \Vert \partial_2^k \theta \Vert_{L^2} +\Vert \partial_2^k u_1 \Vert_{L^{4}_{x_1}(L^2_{x_2})} \Vert \partial_1 \theta \Vert_{L^4_{x_1}(L^{\infty}_{x_2})} \Vert \partial_2^k \theta \Vert_{L^2} \\ 
\triangleq & B_{311} + B_{312}.
\end{align*}
For the first term $B_{311}$, by Minkowski's inequality, interpolation and Young's inequality, we get  
\begin{align*}
B_{311} \lesssim & \Vert \partial_2 u_1 \Vert_{L^4}^{\frac{3}{4}}  \Vert \partial_1\partial_2 u_1 \Vert_{L^4}^{\frac{1}{4}}  \Vert \partial_1 \partial_2^{k-1}\theta \Vert_{L^2}^{\frac{3}{4}}  \Vert \partial_1 \partial_2^{k} \theta  \Vert_{L^2}^{\frac{1}{4}} \Vert \partial_2^k \theta \Vert_{L^2}\\
\lesssim &  \Vert \partial_2 u_1 \Vert_{L^4}^{\frac{3}{2}}  \Vert \partial_1\partial_2 u_1 \Vert_{L^4}^{\frac{1}{2}}  \Vert \partial_2^k \theta \Vert_{L^2}  +  \Vert \partial_1 \partial_2^{k-1} \theta \Vert_{L^2}^{\frac{3}{2}}  \Vert \partial_1 \partial_2^{k}\theta \Vert_{L^2}^{\frac{1}{2}} \Vert \partial_2^k \theta \Vert_{L^2} \\
\lesssim &  \Vert  u_1 \Vert_{L^4}^{\frac{3k-6}{2k-2}} \Vert  \partial_2^{k-1} u_1 \Vert_{L^4}^{\frac{3}{2k-2}}  \Vert \partial_1 u_1 \Vert_{L^4}^{\frac{k-3}{2k-4}}  \Vert \partial_1 \partial_2^{k-2}u_1 \Vert_{L^4}^{\frac{1}{2k-4}}  \Vert \partial_2^k \theta \Vert_{L^2} + \Vert \theta \Vert_{H^k} \Vert \partial_1 \theta \Vert_{H^k}^2 \\
\lesssim & \Vert  u_1 \Vert_{L^4}^{\frac{3k-6}{2k-2}}  \Vert \partial_1 u_1 \Vert_{L^4}^{\frac{k-3}{2k-4}}  \Vert  (u |\theta) \Vert_{H^k}^{\frac{2k^2-2k-3}{2k^2 - 6k +4}} + \Vert \theta \Vert_{H^k} \Vert \partial_1 \theta \Vert_{H^k}^2 .
\end{align*}
For the term $B_{312}$, applying Minkowski's inequality, Young's inequality and related interpolation tools, we can obtain that
\begin{align*}
B_{312} \lesssim & \Vert  \partial_2^k u_1 \Vert_{L^2}^{\frac{3}{4}}  \Vert \partial_1\partial_2^k u_1 \Vert_{L^2}^{\frac{1}{4}}  \Vert \partial_1 \theta \Vert_{L^4}^{\frac{3}{4}}  \Vert \partial_1 \partial_2 \theta\Vert_{L^4}^{\frac{1}{4}} \Vert \partial_2^k \theta \Vert_{L^2} \\
\lesssim & \Vert \partial_2^k u_1 \Vert_{L^2}^{6}  \Vert \partial_1\partial_2^k u_1 \Vert_{L^2}^{2} + \Vert \partial_1 \theta \Vert_{L^4}^{\frac{6}{7}}  \Vert \partial_1 \partial_2 \theta \Vert_{L^4}^{\frac{2}{7}} \Vert \partial_2^k \theta \Vert_{L^2}^{\frac{8}{7}}\\
\lesssim &    \Vert \partial_1 \theta \Vert_{L^4}^{\frac{56k-126}{49k-98}}   \Vert \partial_1 \partial_2^{k-2} \theta \Vert_{L^4}^{\frac{2}{7k-14}} \Vert \partial_2^k \theta \Vert_{L^2}^{\frac{8}{7}}  +  \Vert u \Vert_{H^k}^6 \Vert \partial_1 u \Vert_{H^k}^2\\
\lesssim &  \Vert \partial_1 \theta \Vert_{L^4}^{\frac{56k-126}{49k-98}}   \Vert \theta \Vert_{H^k}^{\frac{56k-98}{49k-98}}   +  \Vert u \Vert_{H^k}^6 \Vert \partial_1 u \Vert_{H^k}^2.
\end{align*}
For $B_{32}$, the estimation follows the same lines as that for $B_{31}$.
\begin{align*}
B_{32}  \lesssim  & \Vert (\partial_2^k (u_2 \partial_2  \theta)- u_2\partial_2^{k+1}\theta) \Vert_{L^2} \Vert \partial_2^k \theta \Vert_{L^2} \\
\lesssim  & \Vert  \partial_1 u_1 \Vert_{L^4}^{\frac{3}{4}}  \Vert \partial_1\partial_2 u_1 \Vert_{L^4}^{\frac{1}{4}}  \Vert \partial_2^{k}\theta \Vert_{L^2}^{\frac{3}{4}}  \Vert \partial_1 \partial_2^{k}\theta  \Vert_{L^2}^{\frac{1}{4}} \Vert \partial_2^k \theta \Vert_{L^2} \\
& +  \Vert \partial_1 \partial_2^{k-1} u_1 \Vert_{L^2}^{\frac{3}{4}}  \Vert \partial_1\partial_2^k u_1 \Vert_{L^2}^{\frac{1}{4}}  \Vert \partial_2 \theta \Vert_{L^4}^{\frac{3}{4}}  \Vert \partial_1 \partial_2 \theta \Vert_{L^4}^{\frac{1}{4}} \Vert \partial_2^k \theta  \Vert_{L^2} \\
\lesssim &   \Vert \partial_1 u_1 \Vert_{L^4}^{\frac{56k-126}{49k-98}}   \Vert (u|\theta)  \Vert_{H^k}^{\frac{56k-98}{49k-98}}   +  \Vert \theta \Vert_{H^k}^6 \Vert \partial_1 \theta \Vert_{H^k}^2 \\
&  + \Vert  \theta \Vert_{L^4}^{\frac{3k-6}{2k-2}} \Vert \partial_1 \theta \Vert_{L^4}^{\frac{k-3}{2k-4}}   \Vert  \theta \Vert_{H^k}^{\frac{2k^2-2k-3}{2k^2 - 6k +4}} + \Vert u \Vert_{H^k} \Vert \partial_1 u \Vert_{H^k}^2.
\end{align*}
Using the same argument as $B_3$, we can bound $B_4$ by:
\begin{align*}
B_4 \lesssim &   \Vert \partial_1 u_1 \Vert_{L^4}^{\frac{56k-126}{49k-98}}   \Vert u  \Vert_{H^k}^{\frac{56k-98}{49k-98}}  +  \Vert u \Vert_{H^k}^6 \Vert \partial_1 u \Vert_{H^k}^2 \\
&  + \Vert  u \Vert_{L^4}^{\frac{3k-6}{2k-2}} \Vert \partial_1 u \Vert_{L^4}^{\frac{k-3}{2k-4}}   \Vert  u \Vert_{H^k}^{\frac{2k^2-2k-3}{2k^2 - 6k +4}}  + \Vert u \Vert_{H^k} \Vert \partial_1 u \Vert_{H^k}^2.
\end{align*}
Thus collecting these above estimates, we obtain the energy inequality:
\begin{equation}\label{energy inequality}
\begin{aligned}
& \Vert (u, \theta)(t) \Vert_{H^k}^2  + 2\int_{0}^{t} \Vert \partial_1 u(\tau) \Vert_{H^k}^2 d \tau  + 2\int_{0}^{t}\Vert \partial_1 \theta(\tau) \Vert_{H^k}^2 d\tau \\
\leq & \Vert (u_0, \theta_0) \Vert_{H^k}^2  + \int_{0}^{t}  C\Vert (u|\theta) \Vert_{H^k} \Vert  (\partial_1u|\partial_1\theta)  \Vert_{H^k}^2 + C\Vert (u|\theta)  \Vert_{H^k}^6 \Vert  (\partial_1u|\partial_1\theta)  \Vert_{H^k}^2 d \tau  \\
& + \int_{0}^{t} C \Vert  (u| \theta) \Vert_{L^4}^{\frac{3k-6}{2k-2}} \Vert (\partial_1 u | \partial_1 \theta) \Vert_{L^4}^{\frac{k-3}{2k-4}}   \Vert  (u | \theta) \Vert_{H^k}^{\frac{2k^2-2k-3}{2k^2 - 6k +4}}  + C\Vert  (\partial_1u_1  | \partial_1\theta) \Vert_{L^4}^{\frac{56k-126}{49k-98}}   \Vert  \theta \Vert_{H^k}^{\frac{56k-98}{49k-98}}d \tau .
\end{aligned}
\end{equation}
By substituting (\ref{equation15})—(\ref{equation16}) into $\eqref{energy inequality}$ and  combining the initial condition (\ref{Initial data condition}), we can get 
\begin{equation}\label{end}
\begin{aligned}
	&\Vert (u| \theta)(t) \Vert_{H^k}^2 + 2\int_{0}^{t} \Vert \partial_1 u(\tau) \Vert_{H^k} ^2 d \tau + 2\int_{0}^{t} \Vert \partial_1 \theta(\tau) \Vert_{H^k} ^2 d \tau \\ \leq  	& C \varepsilon \int_{0}^{t} \Vert \partial_1 u(\tau) \Vert_{H^k} ^2 d \tau+ C \varepsilon \int_{0}^{t} \Vert \partial_1 \theta(\tau) \Vert_{H^k} ^2 d \tau+\varepsilon^2 + 2 CC_0^3 \varepsilon^3 +2CC_0^{\frac{16}{7}} \varepsilon^{\frac{16}{7}}.
	\end{aligned}
\end{equation}
	For $\varepsilon$ sufficiently small such that 
	\begin{equation*}
	C \varepsilon \leq 1, \   2 CC_0^3 \varepsilon \leq 1,\  2CC_0^{\frac{16}{7}} \varepsilon^{\frac{2}{7}} \leq 1,
	\end{equation*}
	we can deduce, taking the universal constant $C_0 \geq 6$,
	\begin{equation*}
	\Vert (u | \theta)(t) \Vert_{H^k}^2 + \int_{0}^{t} \Vert \partial_1 u(\tau) \Vert_{H^k} ^2 d \tau + \int_{0}^{t} \Vert \partial_1 \theta(\tau) \Vert_{H^k} ^2 d \tau  \leq 3 \varepsilon^2 ,
	\end{equation*}
	that is 
	\begin{equation*}
	\Vert (u | \theta)(t) \Vert_{H^k} + \Big(\int_{0}^{t} \Vert \partial_1 u(\tau) \Vert_{H^k} ^2 d \tau \Big)^{\frac{1}{2}} + \Big( \int_{0}^{t} \Vert \partial_1 \theta(\tau) \Vert_{H^k} ^2 d \tau \Big)^{\frac{1}{2}} \leq 3\varepsilon \leq \frac{C_0}{2} \varepsilon,
	\end{equation*}
thus we finish the proof of (\ref{IEu}). In the next section, we will establish the decay estimate (\ref{IDu1})—(\ref{IDpartial2u}).
\section{Decay estimate}\label{section5}
To establish the decay estimate of $u$ and $\theta$, let us recall the formulations of $u$ and $\theta$,
\begin{align*}
u_1 &= \sum_{\sigma \in \{ \pm \}}e^{t\partial_1^2}e^{\sigma t \mathcal{R}_1 } Q^{u_1}_{\sigma} (u_0,\theta_0) + \int_{0}^{t}e^{(t-s)\partial_1^2}e^{\sigma (t-s) \mathcal{R}_1 } Q^{u_1}_{\sigma} (\mathbb{P}(u \cdot \nabla u), u\cdot \nabla \theta ) \ ds , \\
u_2 &= \sum_{\sigma \in \{ \pm \}}e^{t\partial_1^2}e^{\sigma t \mathcal{R}_1 } Q^{u_2}_{\sigma} (u_0,\theta_0) + \int_{0}^{t}e^{(t-s)\partial_1^2}e^{\sigma (t-s) \mathcal{R}_1 } Q^{u_2}_{\sigma} (\mathbb{P}(u \cdot \nabla u), u\cdot \nabla \theta  ) \ ds,\\
\theta &= \sum_{\sigma \in \{ \pm \}}e^{t\partial_1^2}e^{\sigma t \mathcal{R}_1 } Q^{\theta}_{\sigma} (u_0,\theta_0) + \int_{0}^{t}e^{(t-s)\partial_1^2}e^{\sigma (t-s) \mathcal{R}_1 } Q^{\theta}_{\sigma} (\mathbb{P}(u \cdot \nabla u), u\cdot \nabla \theta ) \ ds,
\end{align*}
where $Q^{u_1}$, $Q^{u_2}$ and $Q^{\theta}$ are defined in \eqref{formula of solution4}--\eqref{formula of solution6}.
\subsection{$L^2$ estimate of $\partial_1^{\alpha}u_1$ and $\partial_1^\alpha\theta$}
	Taking the $L^2$ norm of $\partial_1^\alpha u_1$, using the boundedness of the Riesz transform and the $L^2$ isometry of the Fourier transform:
	\begin{align*}
	&
	\Vert \partial_1^\alpha(u_1 | \theta) \Vert_{L^{2}}  \\ \lesssim  & \Vert e^{t\partial_1^2}e^{\sigma t \mathcal{R}_1 } \partial_1^\alpha Q^{u_1}_{\sigma} (u_0,\theta_0) 
	\Vert_{L^{2}} + \int_{0}^{t} \Vert e^{(t-s)\partial_1^2}e^{\sigma (t-s) \mathcal{R}_1 }  \partial_1^\alpha Q^{u_1}_{\sigma} (\mathbb{P}(u \cdot \nabla u), u\cdot \nabla \theta ) \Vert_{L^{2}} d\tau \\
	& +  \Vert e^{t\partial_1^2}e^{\sigma t \mathcal{R}_1 } \partial_1^\alpha Q^{\theta}_{\sigma} (u_0,\theta_0) 
	\Vert_{L^{2}} + \int_{0}^{t} \Vert e^{(t-s)\partial_1^2}e^{\sigma (t-s) \mathcal{R}_1 } \partial_1^\alpha Q^{\theta}_{\sigma} (\mathbb{P}(u \cdot \nabla u), u\cdot \nabla \theta ) \Vert_{L^{2}} d\tau \\
	\lesssim  &  \Vert e^{t\partial_1^2 } \partial_1^\alpha (u_0 | \theta_0) \Vert_{L^{2}} + \int_{0}^{t} \Vert e^{(t-\tau)\partial_1^2 } \partial_1^\alpha (u \cdot \nabla \theta ) \Vert_{L^{2}} d\tau + \int_{0}^{t} \Vert e^{(t-\tau)\partial_1^2 } \partial_1^\alpha (u \cdot \nabla u_1 ) \Vert_{L^{2}} d\tau \\
	& +  \int_{0}^{t} \Vert e^{(t-\tau)\partial_1^2 } \partial_1^\alpha (u \cdot \nabla u_2 ) \Vert_{L^{2}} d\tau\\
	\triangleq & \Vert e^{t\partial_1^2 }\partial_1^\alpha(u_0 | \theta_0) \Vert_{L^{2}} + D_1 + D_2 + D_3.
	\end{align*}
By Lemma \ref{HEAT DECAY LEMMA} and interpolation, the linear term admits a straightforward decay estimate,
	\begin{align*}
	\Vert e^{t\partial_1^2 } \partial_1^\alpha(u_0 | \theta_0)  \Vert_{L^{2}} & \lesssim (1+t)^{-\frac{1+2\alpha}{4} }(\Vert (u_0 | \theta_0) \Vert_{ L^{2}_{x_2} (L^{1}_{x_1}) } + \Vert \partial_1^\alpha(u_0 | \theta_0) \Vert_{L^2} ) \\
	&\lesssim  (1+t)^{-\frac{1+2\alpha}{4} } (\Vert (u_0 | \theta_0) \Vert_{W^{2,1}} + \Vert (u_0 | \theta_0) \Vert_{H^2} ) \\
	& \lesssim  \varepsilon  (1+t)^{-\frac{1+\alpha}{4} }.
	\end{align*}
	For the  term $D_1$, by Lemma \ref{HEAT DECAY LEMMA}, we get
	\begin{align*}
	&\int_{0}^{t} \Vert e^{\partial_1^2 (t-\tau)}  \partial_1^\alpha (u \cdot \nabla \theta) \Vert_{L^{2}} d\tau \\
	\lesssim & \int_{0}^{t} \Vert e^{\partial_1^2 (t-\tau)} \partial_1^\alpha (u_1 \partial_1 \theta+u_2 \partial_2  \theta) \Vert_{L^{2}} d\tau \\
	\lesssim & \int_{0}^{t}(t - \tau )^{-\frac{1+4\alpha}{8}-\eta}( \Vert   u_1 \partial_1 \theta  \Vert_{ L^2_{x_2} (L^{\frac{4}{3+8 \eta}}_{x_1})} +\Vert   u_2 \partial_2 \theta  \Vert_{ L^2_{x_2} (L^{ \frac{4}{3+8 \eta}}_{x_1})} ) d\tau\\
	\triangleq & \int_{0}^{t} (t- \tau)^{-\frac{1+4\alpha}{8} - \eta } (D_{11} + D_{12})(\tau) d\tau.
	\end{align*}
For the term $D_{11}$, by  Minkowski's inequality, H\"older's inequality and interpolation, we can obtain that
	\begin{align*}
	D_{11}(\tau)
	\lesssim & \Big\Vert \Vert u_1 \Vert_{L^{4}_{x_2}} \Vert \partial_1 \theta \Vert_{L^{4}_{x_2}} \Big\Vert_{L^{ \frac{4}{3+8 \eta}}_{x_1}}  \\
	\lesssim &  \Vert u_1 \Vert_{L^{\frac{4}{1+8 \eta}}_{x_1}(L^{4}_{x_2})} \Vert \partial_1 \theta \Vert_{L^{2}_{x_1}(L^{4}_{x_2}) }  \\
	\lesssim & \Vert u_1 \Vert_{L^{\frac{4}{1+8 \eta}}}^{1-2\eta} \Vert \partial_2 u_1 \Vert_{L^{\frac{4}{1+8 \eta}}}^{2\eta} \Vert \partial_1 \theta \Vert_{L^{2}}^{\frac{3}{4}} \Vert \partial_1 \partial_2 \theta \Vert_{L^{2}}^{\frac{1}{4} } \\
	\lesssim & \Vert u_1 \Vert_{L^{4}}^{(1-8\eta)(1-2\eta)}\Vert u_1 \Vert_{L^{2}}^{8\eta (1-2\eta)} \Vert \partial_2 u_1 \Vert_{L^{4}}^{(1-8\eta)2\eta}\\
		&  \times \Vert \partial_2 u_1 \Vert_{L^{2}}^{8\eta \cdot 2\eta}\Vert \partial_1 \theta \Vert_{L^{2}}^{\frac{3}{4}} \Vert \partial_1 \partial_2 \theta \Vert_{L^{2}}^{\frac{1}{4} }\\
	\lesssim & \Vert u_1 \Vert_{L^{4}}^{(1-8\eta)(1-2\eta)}\Vert u_1 \Vert_{L^{2}}^{8\eta (1-2\eta)} \Vert u_1 \Vert_{L^{4}}^{(1-8\eta)2\eta\frac{k-2}{k-1}} \Vert \partial_2^{k-1}u_1 \Vert_{L^{4}}^{(1-8\eta)2\eta\frac{1}{k-1}} \\
	&  \times \Vert \partial_2 u_1 \Vert_{L^{2}}^{8\eta \cdot 2\eta}\Vert \partial_1 \theta \Vert_{L^{2}}^{\frac{3}{4}} \Vert \partial_1 \partial_2 \theta \Vert_{L^{2}}^{\frac{1}{4} }\\
	\lesssim & C_0^2 \varepsilon^2(1+\tau)^{-s_1},
	\end{align*}
	where 
	$$s_1 = \frac{1}{2}(1-8\eta)(1-2\eta) + \frac{k-2}{2k-2}(1-8\eta)2\eta +(\frac{1}{8}+\eta)8\eta+ \frac{5}{8}+ \eta.$$  
For $\eta = \frac{1}{120}$, $k=14$, we can get  $s_1 \thickapprox 1.1088$  which is larger than 1.\\
For $D_{12}$, based on Minkowski's inequality, H\"older's inequality and related interpolation tools, we get
	\begin{align*}
	D_{12}(\tau)  
	\lesssim & \Big\Vert \Vert u_2 \Vert_{L^{\infty}_{x_2}} \Vert \partial_2 \theta \Vert_{L^{2}_{x_2}} \Big\Vert_{L^{\frac{4}{3+8 \eta}}_{x_1}} \\
	\lesssim &  \Vert u_2 \Vert_{L^{\frac{4}{1+8\eta}}_{x_1}(L^{\infty}_{x_2} )} \Vert \partial_2 \theta \Vert_{L^2}  \\
	\lesssim &  \Vert u_2 \Vert_{L^{\frac{4}{1+8\eta}}_{x_1}(L^{\infty}_{x_2} )} \Vert \partial_2 \theta \Vert_{L^2}  \\
	\lesssim &   \Vert u_2 \Vert_{L^{4}}^{(\frac{3}{4} - 2 \eta) (1-8\eta)}  \Vert u_2 \Vert_{L^{2}}^{(\frac{3}{4}- 2 \eta ) 8 \eta} \Vert \partial_1 u_1 \Vert_{L^{4}}^{(\frac{1}{4}+ 2\eta) (1-8\eta)}   \Vert \partial_1 u_1 \Vert_{L^{2}}^{(\frac{1}{4} + 2\eta )  8\eta}\Vert  \partial_2  \theta \Vert_{L^{2}}  \\
	\lesssim &  C_0^2 \varepsilon^2 (1+\tau)^{-s_2},
	\end{align*}
	where $$s_2 = [(\frac{7}{8} + \eta - \delta)(1-8\eta)+ 2\eta](\frac{3}{4} - 2\eta) + [(1-8\eta)(1-\delta)+ 5\eta + 8\eta^2](\frac{1}{4} + 2\eta) +(\frac{1}{8}+ \eta) .$$ 
For $\eta = \frac{1}{120}$, $k=14$ and $\delta = 10^{-5}$, we can get  $s_2 \thickapprox 1.01028$  which is larger than 1.\\
	Therefore by Lemma \ref{decay lem}, we derive that 
	\begin{align*}
	D_1 \lesssim &  C_0^2 \varepsilon^2 \int_{0}^{t} (t - \tau )^{-\frac{1+4\alpha}{8}-\eta} \Big((1+\tau)^{- s_1}+(1+\tau)^{- s_2 } \Big) d \tau \\
	\lesssim &  C_0^2 \varepsilon^2 (1+t)^{-\frac{1+4\alpha}{8}-\eta} .
	\end{align*}
	The   term $D_2$ can be estimated in a similar way as for $D_1$ that,
	\begin{equation*}
	D_2 \lesssim  C_0^2 \varepsilon^2 (1+t)^{-\frac{1+4\alpha}{8}-\eta}.
	\end{equation*}
	For the nonlinear term $D_3$, since the decay of $u_2$ is faster than that of $u_1$ in both the $L^2$ and $L^4$ norms, the estimate of $D_3$ is easier than that of $D_1$; we omit the details.
	\begin{equation*}
		D_3 \lesssim  C_0^2 \varepsilon^2 (1+t)^{-\frac{1+4\alpha}{8}-\eta}.
	\end{equation*}
	Combining the estimates above, we obtain
	\begin{align*}
	\Vert (u_1 | \theta) \Vert_{L^{2}} \lesssim(1+C_0^2 \varepsilon) \varepsilon (1+t)^{-\frac{1+4\alpha}{8}- \eta}.
	\end{align*}
Taking $\varepsilon$ sufficiently small, one can derive the bootstrap
conclusion (\ref{IDu1}).
	\subsection{$L^4$ estimate of $\partial_1^\alpha u_1$ and $\partial_1^\alpha \theta$}
  Taking the $L^4$ norm  of $\partial_1^\alpha u_1$, using the $L^4$ boundedness of the Riesz transform, we  get
	\begin{align*}
	&\Vert \partial_1^\alpha (u_1 | \theta) \Vert_{L^{4}}   \\\lesssim& \Vert e^{t\partial_1^2}e^{\sigma t \mathcal{R}_1 } \partial_1^\alpha  Q^{u_1}_{\sigma} (u_0,\theta_0) 
	\Vert_{L^{4}} + \int_{0}^{t} \Vert e^{(t-\tau)\partial_1^2}e^{\sigma (t-\tau) \mathcal{R}_1 } \partial_1^\alpha  Q^{u_1}_{\sigma} (\mathbb{P}(u \cdot \nabla u), u\cdot \nabla \theta  ) \Vert_{L^{4}} d\tau \\
	& + \Vert e^{t\partial_1^2}e^{\sigma t \mathcal{R}_1 } \partial_1^\alpha  Q^{\theta}_{\sigma} (u_0,\theta_0) 
	\Vert_{L^{4}} + \int_{0}^{t} \Vert e^{(t-\tau)\partial_1^2}e^{\sigma (t-\tau) \mathcal{R}_1 } \partial_1^\alpha  Q^{\theta}_{\sigma} (\mathbb{P}(u \cdot \nabla u), u\cdot \nabla \theta  ) \Vert_{L^{4}} d\tau \\
	\lesssim   & \Vert e^{t\partial_1^2 }e^{\sigma t \mathcal{R}_1} \partial_1^\alpha (u_0 | \theta_0) \Vert_{L^{4}} + \int_{0}^{t} \Vert e^{(t-\tau)\partial_1^2 } e^{\sigma (t-\tau) \mathcal{R}_1 } \partial_1^\alpha  (u \cdot \nabla u | u\cdot \nabla \theta ) \Vert_{L^{4}} d\tau.
	\end{align*}
By Proposition \ref{3}, Lemma \ref{Localized heat kernel}
 and interpolation, the linear term can be estimated:
	\begin{align*}
	\Vert e^{t\partial_1^2 }e^{\sigma t \mathcal{R}_1 }\partial_1^\alpha  (u_0 | \theta_0) \Vert_{L^{4}} &\lesssim  (1+t)^{-\frac{1}{2}}  \Vert \partial_1^\alpha e^{\frac{t\partial_1^2}{2} }(u_0 | \theta_0) \Vert_{W^{1+\gamma, \frac{4}{3}}} \\
	&\lesssim (1+t)^{-\frac{1+\alpha}{2}} (\Vert (u_0 | \theta_0) \Vert_{W^{1+\gamma, \frac{4}{3}}}+ \Vert \partial_1^\alpha(u_0 | \theta_0) \Vert_{W^{1+\gamma, \frac{4}{3}}})\\
	& \lesssim (1+t)^{-\frac{1+\alpha}{2}} (\Vert (u_0 | \theta_0) \Vert_{W^{3,1}} + \Vert (u_0 | \theta_0) \Vert_{H^3} ) \\
	& \lesssim  \varepsilon  (1+t)^{-\frac{1+\alpha}{2}}.
	\end{align*}	%
	For the nonlinear terms, by Proposition \ref{3}, Lemma \ref{Localized heat kernel} and H\"older's inequality, we get 
	\begin{align*}
	&\int_{0}^{t} \Vert e^{(t-\tau)\partial_1^2 } e^{\sigma (t-\tau) \mathcal{R}_1 } \partial_1^\alpha (u \cdot \nabla u | u\cdot \nabla \theta ) \Vert_{L^{4}} d\tau \\
	\lesssim & \int_{0}^{t} (t- \tau)^{-\frac{1+(1-\delta)\alpha}{2}} \Vert  (u \cdot \nabla u | u\cdot \nabla \theta ) \Vert_{W^{1+2\gamma,\frac{4}{3}}}  d\tau \\
	\lesssim & \int_{0}^{t} (t- \tau)^{-\frac{1+(1-\delta)\alpha}{2} } \bigg(\Vert  (u \cdot \nabla u | u\cdot \nabla \theta ) \Vert_{L^{\frac{4}{3}}} + \Vert  \Lambda_2^{1+2\gamma}(u \cdot \nabla u | u\cdot \nabla \theta ) \Vert_{L^{\frac{4}{3}}}\\&\quad+\Vert  \Lambda_1^{1+2\gamma}(u \cdot \nabla u | u\cdot \nabla \theta) \Vert_{L^{\frac{4}{3}}} d\tau \bigg)\\ 
	\triangleq & \int_{0}^{t}(t-\tau)^{-\frac{1+(1-\delta)\alpha}{2}} (E_1 + E_2 + E_3) d\tau.
	\end{align*} 
	For the  term $E_1$, by  H\"older's inequality, we have
	\begin{align*}
	E_1(\tau) & \lesssim \Vert u_1 \Vert_{L^{4}}  \Vert (\partial_1 u | \partial_1 \theta) \Vert_{L^{2}} + \Vert u_2 \Vert_{L^{4}}  \Vert (\partial_2 u | \partial_2 \theta) \Vert_{L^{2}} \\
	& \lesssim C_0^2 \varepsilon^2 (1+\tau)^{-(1 +2\eta - \delta )}.
	\end{align*}
We now turn to the term $E_2$ and decompose it into three parts.
\begin{align*}
	E_2(\tau)  = & 	\Vert \Lambda_2^{1+2\gamma}(u \cdot \nabla \theta) \Vert_{L^{\frac{4}{3}}} 
	+ \Vert \Lambda_2^{1+2\gamma}(u \cdot \nabla u_1) \Vert_{L^{\frac{4}{3}}} + 	\Vert \Lambda_2^{1+2\gamma}(u \cdot \nabla u_2) \Vert_{L^{\frac{4}{3}}} \\
	\triangleq & E_{21} + E_{22} + E_{23}.
\end{align*}
By Lemma \ref{Kato-Ponce inequality} and H\"older's inequality, we can get
	\begin{align*}
	E_{21} \lesssim &\Vert u_1 \Vert_{L^{4}}  \Vert \Lambda_2^{1+2\gamma} \partial_1 \theta \Vert_{L^{2}} + \Vert \Lambda_2^{1+2\gamma}u_1 \Vert_{L^{2}}  \Vert \partial_1 \theta \Vert_{L^{4}} \\
	& + \Vert u_2 \Vert_{L^{4}}  \Vert \Lambda_2^{1+2\gamma} \partial_2 \theta \Vert_{L^{2}} + \Vert \Lambda_2^{1+2\gamma}  u_2 \Vert_{L^{2}}  \Vert \partial_2 \theta \Vert_{L^{4}}.
	\end{align*}
	By interpolation,
	\begin{align*}
	\Vert \Lambda_2^{2+2\gamma} \theta \Vert_{L^{2}} \lesssim \Vert  \partial_2 \theta \Vert_{L^{2}}^{\frac{k-2-2\gamma}{k-1}} \Vert \partial_2^k \theta \Vert_{L^{2}}^{\frac{1+2\gamma}{k-1}}.
	\end{align*}
	The other terms can be estimated through interpolation similarly. 
	Since $k \geq 14 $ and $\gamma= 10^{-5}$, by the assumption of decay estimate (\ref{equation15})—(\ref{equation16}), we get
	\begin{equation*}
	E_{21} \lesssim C_0^2 \varepsilon^2 (1+\tau)^{-1- \delta } .
	\end{equation*}
	The term $E_{22}$ contained in $E_2$ can be estimated similarly.\\
	\begin{equation*}
		E_{22} \lesssim C_0^2 \varepsilon^2 (1+\tau)^{-1- \delta } .
	\end{equation*}
	For the nonlinear term $\Vert \Lambda_2^{1+2\gamma}(u \cdot \nabla u_2) \Vert_{L^{\frac{4}{3}}}$ contained in $E_2$, by  H\"older's inequality  and Lemma \ref{Kato-Ponce inequality}, we get
	\begin{align*}
	\Vert \Lambda_2^{1+2\gamma}(u \cdot \nabla u_2) \Vert_{L^{\frac{4}{3}}} \lesssim &\Vert u_1 \Vert_{L^{4}}  \Vert \Lambda_2^{1+2\gamma} \partial_1 u_2 \Vert_{L^{2}} + \Vert \Lambda_2^{1+2\gamma} u_1 \Vert_{L^{2}}  \Vert \partial_1 u_2 \Vert_{L^{4}} \\
	& + \Vert u_2 \Vert_{L^{4}}  \Vert \Lambda_2^{1+2\gamma}\partial_2 u_2 \Vert_{L^{2}} + \Vert  \Lambda_2^{1+2\gamma} u_2 \Vert_{L^{2}}  \Vert \partial_2 u_2 \Vert_{L^{4}}.
	\end{align*}
	Then by interpolation, we can get 
	\begin{equation*}
	E_{23} \lesssim C_0^2 \varepsilon^2 (1+\tau)^{-1-\delta}.
	\end{equation*}
	For $E_3$, because the dissipation is in the horizontal direction, the estimate is easier compared with $E_2$. We can obtain that
	\begin{equation*}
		E_3(\tau) \lesssim C_0^2 \varepsilon^2 (1+t)^{-1-\delta}.
	\end{equation*}
 Combining the estimates above, by Lemma \ref{decay lem}, we obtain
	\begin{align*}
	\Vert (u_1 | \theta) \Vert_{L^{4}} \lesssim(1+C_0^2 \varepsilon) \varepsilon (1+t)^{-\frac{1+(1-\delta)\alpha}{2} }.
	\end{align*}
Taking $\varepsilon$ sufficiently small, one can derive the bootstrap
conclusion (\ref{IDuL4}).	
	\subsection{$L^2$ estimate of $\partial_1^\alpha u_2$}
\noindent Due to the special structure of the nonlinear equation, both the $L^2$ norm and the $L^4$ norm of $u_2$ enjoy a faster decay rate than $u_1$. Taking the $L^2$ norm to  $\partial_1^\alpha u_2$, using the boundedness of the Riesz transform and the $L^2$ isometry of the Fourier transform:
	\begin{align*}
	\Vert \partial_1^\alpha u_2 \Vert_{L^{2}} & \lesssim \Vert e^{t\partial_1^2}e^{\sigma t \mathcal{R}_1 } \partial_1^\alpha  Q^{u_2}_{\sigma} (u_0,\theta_0)\Vert_{L^{2}} + \int_{0}^{t} \Vert e^{(t-\tau)\partial_1^2}e^{\sigma (t-\tau) \mathcal{R}_1 } \partial_1^\alpha Q^{u_2}_{\sigma} (\mathbb{P}(u \cdot \nabla u), u \cdot \nabla \theta ) \Vert_{L^{2}} d\tau \\
	& \lesssim  \Vert e^{t\partial_1^2 } \partial_1^\alpha (u_0 | \theta_0) \Vert_{L^{2}} +\int_{0}^{t} \Vert e^{(t-\tau)\partial_1^2 } \partial_1^\alpha \partial_1 (u \theta )   \Vert_{L^{2}} d\tau + \int_{0}^{t} \Vert e^{(t-\tau)\partial_1^2 } \partial_1^\alpha \partial_1 (u \otimes u )   \Vert_{L^{2}} d\tau \\
	& \triangleq \Vert e^{t\partial_1^2 } \partial_1^\alpha (u_0 | \theta_0) \Vert_{L^{2}} + F_1 + F_2.
	\end{align*}
	The linear term admits a straightforward decay estimate:
	\begin{align*}
	\Vert e^{t \partial_1^2 } \partial_1^\alpha (u_0 | \theta_0) \Vert_{L^{2}} &\lesssim (1+t)^{-\frac{1+2\alpha}{4} }(\Vert  (u_0 | \theta_0) \Vert_{L^{1}_{x_1} L^{2}_{x_2}} + \Vert \partial_1^\alpha (u_0 | \theta_0) \Vert_{L^2} ) \\
	& \lesssim (1+t)^{-\frac{1+2\alpha}{4}} (\Vert (u_0 | \theta_0) \Vert_{W^{2,1}} + \Vert (u_0 | \theta_0) \Vert_{H^2} ) \\
	& \lesssim  \varepsilon  (1+t)^{-\frac{1+2\alpha}{4} }.
	\end{align*}
	For the nonlinear terms $F_1$,
	\begin{align*}
F_1 & \lesssim \int_{0}^{t} \Vert e^{\partial_1^2 (t-\tau)} \partial_1^\alpha  \partial_1(u_1 \theta)  \Vert_{L^{2}} d\tau + \int_{0}^{t} \Vert e^{\partial_1^2 (t-\tau)} \partial_1^\alpha  \partial_1 (u_2 \theta) \Vert_{L^{2}} d\tau\\
  & \triangleq F_{11} + F_{12}.
	\end{align*}
	By the boundedness of the Hilbert transform and Lemma \ref{HEAT DECAY LEMMA}, we have
	\begin{align*}
	 F_{11}
	\lesssim & \int_{0}^{t} \Vert e^{\partial_1^2 (t-\tau)}  \Lambda_1^{1+\alpha}(u_1  \theta + u_2 \theta) \Vert_{L^{2}} d\tau \\
	\lesssim & \int_{0}^{t}(t - \tau )^{-\frac{1+2\alpha}{4}} (\Vert   \Lambda_1^{\frac{3}{4}}(u_1  \theta) \Vert_{L^2_{x_2} (L^{\frac{4}{3}}_{x_1} ) }+ \Vert   \Lambda_1^{\frac{3}{4}}(u_2  \theta) \Vert_{L^2_{x_2} (L^{\frac{4}{3}}_{x_1} ) }) d\tau \\
	\triangleq & \int_{0}^{t}(t - \tau )^{-\frac{1+2\alpha}{4}} (F_{111} + F_{112})(\tau) d\tau. 
	\end{align*}
Applying Minkowski's inequality, H\"older's inequality together with Lemma \ref{Kato-Ponce inequality}, we derive
	\begin{align*}
	  F_{111}
	\lesssim & \Big\Vert \Vert \Lambda_1^{\frac{3}{4} } \theta \Vert_{L^{4}_{x_1}} \Vert u_1 \Vert_{L^{2}_{x_1}} + \Vert \Lambda_1^{\frac{3}{4} } u_1 \Vert_{L^{4}_{x_1}} \Vert \theta \Vert_{L^{2}_{x_1}} \Big\Vert_{L^{2}_{x_2}}  \\
	\lesssim & \Big\Vert \Vert \Lambda_1 \theta \Vert_{L^{4}_{x_1}}^{\frac{3}{4}   }  \Vert \theta \Vert_{L^{4}_{x_1}}^{\frac{1}{4}   } \Vert u_1 \Vert_{L^{2}_{x_1}} \Big\Vert_{L^{2}_{x_2}}  + \Big\Vert \Vert \Lambda_1 u_1 \Vert_{L^{4}_{x_1}}^{\frac{3}{4}   }  \Vert u_1 \Vert_{L^{4}_{x_1}}^{\frac{1}{4}   } \Vert \theta \Vert_{L^{2}_{x_1}} \Big\Vert_{L^{2}_{x_2}}   \\
	\lesssim &  \Vert \Lambda_1 \theta \Vert_{L^4}^{\frac{3}{4}   }  \Vert \theta \Vert_{L^4}^{\frac{1}{4}  } \Vert u_1 \Vert_{L^{2}_{x_1}(L^4_{x_2}) } + \Vert \Lambda_1 u_1 \Vert_{L^4}^{\frac{3}{4}   }  \Vert u_1 \Vert_{L^4}^{\frac{1}{4}  } \Vert \theta \Vert_{L^{2}_{x_1} (L^4_{x_2}) }.
	\end{align*}
Via interpolation in the $x_2$ variable, we have
	\begin{align*}
	\Vert  \theta \Vert_{L^2_{x_1}(L^4_{x_2})} \leq \Vert  \theta  \Vert_{L^2}^{\frac{3}{4}} \Vert \partial_2 \theta \Vert_{L^2}^{\frac{1}{4}},
	\end{align*}
and
\begin{align*}
\Vert u_1 \Vert_{L^2_{x_1}(L^4_{x_2})} &\leq \Vert u_1 \Vert_{L^2}^{\frac{3}{4}   }  	\Vert \partial_2 u_1 \Vert_{L^2}^{\frac{1}{4} } .
\end{align*}
	Similarly, by Minkowski's inequality, H\"older's inequality and interpolation,
	\begin{align*}
	F_{112}	\lesssim & \Big\Vert \Vert \Lambda_1^{\frac{3}{4} } u_2 \Vert_{L^{4}_{x_1}} \Vert \theta \Vert_{L^{2}_{x_1}} \Big\Vert_{L^{2}_{x_2}}  + \Big\Vert \Vert \Lambda_1^{\frac{3}{4} } \theta \Vert_{L^{4}_{x_1}} \Vert u_2 \Vert_{L^{2}_{x_1}} \Big\Vert_{L^{2}_{x_2}} \\
	\lesssim & \Big\Vert \Vert \Lambda_1 u_2 \Vert_{L^{4}_{x_1}}^{\frac{3}{4}   }  \Vert u_2 \Vert_{L^{4}_{x_1}}^{\frac{1}{4}   }  \Vert \theta \Vert_{L^{2}_{x_1}} \Big\Vert_{L^{2}_{x_2}}  + \Big\Vert \Vert \Lambda_1 \theta \Vert_{L^{4}_{x_1}}^{\frac{3}{4}   }  \Vert \theta \Vert_{L^{4}_{x_1}}^{\frac{1}{4}   }  \Vert u_2 \Vert_{L^{2}_{x_1}} \Big\Vert_{L^{2}_{x_2}} \\
	\lesssim &  \Vert \Lambda_1 u_2 \Vert_{L^4}^{\frac{3}{4}   }  \Vert u_2 \Vert_{L^4}^{\frac{1}{4}  } \Vert \theta \Vert_{L^{2}_{x_1} (L^4_{x_2}) } + \Vert \Lambda_1 \theta \Vert_{L^4}^{\frac{3}{4}   }  \Vert \theta \Vert_{L^4}^{\frac{1}{4}  } \Vert u_2 \Vert_{L^{2}_{x_1} (L^4_{x_2}) }.
	\end{align*}
And	by interpolation in $x_2$,
	\begin{align*}
	\Vert  u_2 \Vert_{L^2_{x_1}(L^4_{x_2})} &\lesssim \Vert  u_2 \Vert_{L^2}^{\frac{3}{4}} \Vert  \partial_2 u_2 \Vert_{L^2}^{\frac{1}{4}}, 
	\end{align*}
 according to the assumption of decay estimate (\ref{equation15})—(\ref{equation16}), we can get
 \begin{equation*}
 F_{111} + F_{112} \lesssim C_0^2 \varepsilon^2 (1+\tau)^{-1-\eta + \delta}.
 \end{equation*}
	Thus by Lemma \ref{decay lem}, we have
	\begin{align*}
F_1	\lesssim & C_0^2 \varepsilon^2 \int_{0}^{t} (t - \tau )^{-\frac{1+2\alpha}{4}} (1+ \tau)^{-1- \eta + \delta } d \tau \\
	\lesssim &C_0^2 \varepsilon^2(1+t)^{-\frac{1+2\alpha}{4}}.
	\end{align*}
	The  nonlinear term $F_2$ can be estimated similarly.
	\begin{equation*}
		F_2	\lesssim C_0^2 \varepsilon^2(1+t)^{-\frac{1+2\alpha}{4}}.
	\end{equation*}
	Combining the estimates above, we obtain
	\begin{align*}
	\Vert u_2 \Vert_{L^{2}} \lesssim(1+C_0^2 \varepsilon) \varepsilon (1+t)^{-\frac{1+2\alpha}{4}}.
	\end{align*}
	Taking $\varepsilon$ sufficiently small, one can derive the bootstrap
	conclusion (\ref{IDu2}).	

	\subsection{$L^4$ estimate of $u_2$}
 Taking the $L^4$ norm to $u_2$, using the $L^4$ boundedness of the Riesz transform, we can get
	\begin{align*}
	\Vert u_2 \Vert_{L^{4}} & \lesssim \Vert e^{t\partial_1^2}e^{\sigma t  \mathcal{R}_1} Q^{u_2}_{\sigma} (u_0,\theta_0)\Vert_{L^{4}} +  \int_{0}^{t} \Vert e^{(t-\tau)\partial_1^2}e^{\sigma (t-\tau) \mathcal{R}_1 } Q^{u_2}_{\sigma} (\mathbb{P}(u \cdot \nabla u), u \cdot \nabla \theta ) \Vert_{L^{4}}  d\tau \\
	& \lesssim \Vert e^{t\partial_1^2}e^{\sigma t  \mathcal{R}_1} \partial_1 \Lambda^{-1} (u_0 | \theta_0)\Vert_{L^{4}} +  \int_{0}^{t} \Vert e^{(t-\tau)\partial_1^2}e^{\sigma (t-\tau) \mathcal{R}_1 } \partial_1 \Lambda^{-1}(u \cdot \nabla u | u \cdot \nabla \theta ) \Vert_{L^{4}}  d\tau \\
	& \triangleq G_1 + G_2.
	\end{align*}
By Proposition \ref{3}, (\ref{localized L4 to L43}) and the divergence-free condition of $u$, the linear term $G_1$ admits the following decay estimate:
\begin{equation}
	\begin{aligned} 
G_1 \lesssim & \sum_{j \in \mathbb{Z}} \Vert  e^{t\partial_1^2}e^{\sigma t \mathcal{R}_1 } \partial_1 \Lambda^{-1}  \Delta_j(u_0| \theta_0) \Vert_{L^{4}} \\
\lesssim & \sum_{j \in \mathbb{Z}} 2^{j} (1+t)^{-\frac{1}{2}} \Vert  e^{t\partial_1^2} \partial_1 \Lambda^{-1}  \Delta_j(u_0| \theta_0) \Vert_{L^{\frac{4}{3}}} \\
& + \sum_{j \in \mathbb{Z}} (1+t)^{-\frac{1}{4}} e^{-c2^{2j}t} 2^{j}\Vert   \Delta_j(u_0| \theta_0) \Vert_{L^{\frac{4}{3}}} \\
\triangleq & G_{11} + G_{12}.
\end{aligned}
\end{equation}
	For the term $G_{11}$, by Lemma \ref{Bernstein} and Lemma \ref{Localized heat kernel}, we have
	\begin{align*} 
	G_{11} = & \sum_{j \in \mathbb{Z}} 2^{j} (1+t)^{-\frac{1}{2}} \Vert  e^{t\partial_1^2} \partial_1 \Lambda^{-1}  \Delta_j(u_0| \theta_0) \Vert_{L^{\frac{4}{3}}} \\
	\lesssim & \sum_{j \in \mathbb{Z}}  (1+t)^{-\frac{1}{2}} \Vert  e^{t\partial_1^2} \partial_1   \Delta_j(u_0| \theta_0) \Vert_{L^{\frac{4}{3}}} \\
	\lesssim & \sum_{j \in \mathbb{Z}}  (1+t)^{-1+\delta} (2^{2j\delta}+2^j) \Vert   \Delta_j(u_0| \theta_0) \Vert_{L^{\frac{4}{3}}} \\
	\lesssim & (1+t)^{-1+\delta}  (\Vert   (u_0| \theta_0) \Vert_{W^{2,1}} + \Vert   (u_0| \theta_0) \Vert_{W^{2,2}} ) .
	\end{align*}
	For the term $G_{12}$, we have 
	\begin{align*} 
	G_{12} 
	\lesssim & \sum_{j \in \mathbb{Z}} (1+t)^{-\frac{1}{4}} e^{-c2^{2j}t} 2^{\frac{3j}{2}}\Vert   \Delta_j(u_0| \theta_0) \Vert_{L^{1}} \\
	\lesssim &  \sum_{j \in \mathbb{Z}} (1+t)^{-1 + \delta}  (2^{\frac{3j}{2}} + 2^{2j\delta})\Vert   \Delta_j(u_0| \theta_0) \Vert_{L^{1}} \\
	\lesssim & (1+t)^{-1+\delta}  (\Vert   (u_0| \theta_0) \Vert_{W^{2,1}} + \Vert   (u_0| \theta_0) \Vert_{W^{2,2}} ) .
	\end{align*}
So collecting the estimate for $G_{11}$ and $G_{12}$, we can get 
\begin{equation}\label{G_1}
G_1 \lesssim \varepsilon (1+t)^{-1+\delta}.
\end{equation}
	For the nonlinear term $G_2$, by Lemma \ref{Bernstein}, we obtain that
	\begin{align*}
	G_2	\lesssim &   \sum_{j \in \mathbb{Z}} \Vert   e^{(t-\tau)\partial_1^2}e^{\sigma (t-\tau) \mathcal{R}_1 } \partial_1 \Lambda_1^{-1} \Delta_j  (u \cdot \nabla u | u \cdot \nabla \theta ) \Vert_{L^{4}}  \\
	\lesssim	&   \sum_{j \in \mathbb{Z}}  (t-\tau)^{-\frac{1}{4}}e^{-c2^{2j}(t - \tau )}2^{j} \Vert  \Delta_j \partial_1 \Lambda^{-1}(u \cdot \nabla u | u \cdot \nabla \theta ) \Vert_{L^{\frac{4}{3}}}  \\
	& +  \sum_{j \in \mathbb{Z}}  (t-\tau)^{-\frac{1}{2}} 2^{j} \Vert  e^{(t-\tau)\partial_1^2}  \partial_1 \Lambda^{-1} \Delta_j(u \cdot \nabla u | u \cdot \nabla \theta ) \Vert_{L^{\frac{4}{3}}} \\
	\triangleq & G_{21} + G_{22}.
	\end{align*}
	For $G_{21}$, by Lemma \ref{Bernstein}, we further deduce from the above that
	\begin{align*}
	G_{21} \lesssim &  \sum_{j \in \mathbb{Z}}  (t-\tau)^{-\frac{1}{4}}e^{-c2^{2j}(t - \tau )}2^{j} \Vert  \Delta_j \partial_1 \Lambda^{-1}(u \cdot \nabla u | u \cdot \nabla \theta ) \Vert_{L^{\frac{4}{3}}} \\
	\lesssim & \sum_{j \in \mathbb{Z}}  (t-\tau)^{-\frac{1}{4}}e^{-c2^{2j}(t - \tau )}2^{j} \Vert  \Delta_j \Lambda_1 (u \otimes u + u \theta ) \Vert_{L^{\frac{4}{3}}} \\
	\lesssim & \sum_{j \in \mathbb{Z}}  (t-\tau)^{-\frac{1}{4}}e^{-c2^{2j}(t - \tau )}2^{\frac{5j}{4}}2^{2\eta j} 2^{-2\delta j}\Vert  \Delta_j \Lambda_1^{\frac{3}{4}- 2\eta + 2\delta} (u \otimes u + u \theta ) \Vert_{L^{\frac{4}{3}}} \\
	\lesssim & (t-\tau)^{-\frac{7}{8} - \eta + \delta} \Vert \Lambda_1^{\frac{3}{4}- 2\eta + 2\delta} (u \otimes u + u \theta ) \Vert_{L^{\frac{4}{3}}} .
	\end{align*}
	By applying  Lemma \ref{Kato-Ponce inequality}, H\"older's inequality and interpolation,
	\begin{align*}
	& \Vert   \Lambda_1^{\frac{3}{4} - 2\eta + 2\delta}(u_1  \theta) \Vert_{L^{\frac{4}{3}} } \\
	\lesssim &  \Vert \Lambda_1^{\frac{3}{4} - 2\eta + 2\delta} \theta \Vert_{L^{2}} \Vert u_1 \Vert_{L^{4}} + \Vert \Lambda_1^{\frac{3}{4} - 2\eta + 2\delta} u_1 \Vert_{L^{2}} \Vert \theta \Vert_{L^{4}}  \\
	\lesssim &  \Vert \Lambda_1 \theta \Vert_{L^2}^{ \frac{3}{4} - 2\eta + 2\delta  }  \Vert \theta \Vert_{L^2}^{\frac{1}{4}  + 2\eta - 2\delta} \Vert u_1 \Vert_{L^{4} } + \Vert \Lambda_1 u_1 \Vert_{L^2}^{\frac{3}{4} - 2\eta + 2\delta  }  \Vert u_1 \Vert_{L^2}^{\frac{1}{4} + 2\eta - 2\delta } \Vert \theta \Vert_{L^{4} } .
	\end{align*}
	Similarly, via Lemma \ref{Kato-Ponce inequality}, H\"older's inequality and standard interpolation procedures, we can get
	\begin{align*}
& \Vert   \Lambda_1^{\frac{3}{4}- 2\eta + 2\delta}(u_2  \theta) \Vert_{L^{\frac{4}{3}} } \\
\lesssim &  \Vert \Lambda_1^{\frac{3}{4} - 2\eta + 2\delta} \theta \Vert_{L^{2}} \Vert u_2 \Vert_{L^{4}} + \Vert \Lambda_1^{\frac{3}{4}- 2\eta + 2\delta } u_2 \Vert_{L^{2}} \Vert \theta \Vert_{L^{4}}  \\
\lesssim &  \Vert \Lambda_1 \theta \Vert_{L^2}^{\frac{3}{4} - 2\eta + 2\delta  }  \Vert \theta \Vert_{L^2}^{\frac{1}{4} + 2\eta - 2\delta } \Vert u_2 \Vert_{L^{4} } + \Vert \Lambda_1 u_2 \Vert_{L^2}^{\frac{3}{4}  - 2\eta + 2\delta }  \Vert u_2 \Vert_{L^2}^{\frac{1}{4} + 2\eta - 2\delta} \Vert \theta \Vert_{L^{4} }.
\end{align*}
 According to the assumption of decay estimate (\ref{equation15})—(\ref{equation16}), we arrive
 \begin{equation*}
 	\Vert \Lambda_1^{\frac{3}{4}- 2\eta + 2\delta} (u \otimes u + u \theta ) \Vert_{L^{\frac{4}{3}}} \lesssim C_0^2 \varepsilon^2(1+\tau)^{-1-\delta}.
 \end{equation*}
	For $G_{22}$, it can be estimated as
\begin{align*}
G_{22} \lesssim &  \sum_{j \in \mathbb{Z}}  (t - \tau )^{-\frac{1}{2}} 2^{j}\Vert e^{(t-\tau)\partial_1^2} \partial_1 \Lambda^{-1} \Delta_j (u \cdot \nabla u | u \cdot \nabla \theta )  \Vert_{L^{\frac{4}{3}}} \\
\lesssim &  \sum_{j \in \mathbb{Z}}  (t - \tau )^{-1+ \delta} 2^{2\delta j}\Vert   \Delta_j (u \cdot \nabla u | u \cdot \nabla \theta )  \Vert_{L^{\frac{4}{3}}} \\
\lesssim & (t- \tau)^{-1+ \delta} (\Vert   (u \cdot \nabla u | u \cdot \nabla \theta )   \Vert_{L^{\frac{4}{3}}} + \Vert   \Lambda_1^{\gamma} (u \cdot \nabla u | u \cdot \nabla \theta )   \Vert_{L^{\frac{4}{3}}}+ \Vert   \Lambda_2^{\gamma} (u \cdot \nabla u |  u \cdot \nabla \theta )  \Vert_{L^{\frac{4}{3}}}).
\end{align*}
By the estimate for $E_1$, $E_2$ and  $E_3$, we have
\begin{equation*}
	\Vert   (u \cdot \nabla u | u \cdot \nabla \theta )   \Vert_{L^{\frac{4}{3}}} + \Vert   \Lambda_1^{\gamma} (u \cdot \nabla u | u \cdot \nabla \theta )   \Vert_{L^{\frac{4}{3}}}+ \Vert   \Lambda_2^{\gamma} (u \cdot \nabla u | u \cdot \nabla \theta )  \Vert_{L^{\frac{4}{3}}} \lesssim C_0^2 \varepsilon^2 (1+\tau)^{-1-\delta}.
\end{equation*}
	Thus, in the end, by Lemma \ref{decay lem}, we get
	\begin{equation}\label{G_2}
	\begin{aligned}
	G_2 \lesssim & C_0^2 \varepsilon^2\int_{0}^{t} (t - \tau)^{-\frac{7}{8} - \eta + \delta} (1+\tau)^{-1- \delta} d \tau \\
	\lesssim & C_0^2 \varepsilon^2 (1+t)^{-\frac{7}{8} - \eta + \delta}.
	\end{aligned}
	\end{equation}
	Combining the estimates above \eqref{G_1}, \eqref{G_2}, we obtain
	\begin{align*}
	\Vert u_2 \Vert_{L^{4}} \lesssim(1+C_0^2 \varepsilon) \varepsilon (1+t)^{- \frac{7}{8}- \eta + \delta}.
	\end{align*}
	Taking $\varepsilon$ sufficiently small, one can derive the bootstrap
	conclusion (\ref{IDu2L4}).	
	\subsection{$L^2$ estimate of $\partial_1^\alpha \partial_2u_1$ and $\partial_1^\alpha \partial_2\theta$}
In this subsection, we will give the $L^2$ decay estimate of $\partial_1^\alpha \partial_2u_1$ and $\partial_1^\alpha \partial_2 \theta$. From the equation, we can get:
	\begin{align*}
	\Vert \partial_1^\alpha (\partial_2u_1 | \partial_2\theta) \Vert_{L^{2}}  \lesssim  & \Vert e^{t\partial_1^2}e^{\sigma t \mathcal{R}_1 } \partial_1^\alpha  Q^{u_1}_{\sigma} (\partial_2u_0,\partial_2\theta_0) 
	\Vert_{L^{2}} +  \Vert e^{t\partial_1^2}e^{\sigma t \mathcal{R}_1} \partial_1^\alpha  Q^{\theta}_{\sigma} (\partial_2u_0,\partial_2\theta_0) \Vert_{L^{2}} \\
	&+ \int_{0}^{t} \Vert e^{(t-\tau)\partial_1^2}e^{\sigma (t-\tau) \mathcal{R}_1} \partial_1^\alpha  Q^{u_1}_{\sigma}  (\partial_2\mathbb{P}(u \cdot \nabla u), \partial_2(u\cdot \nabla \theta)) \Vert_{L^{2}} d\tau \\
	&  + \int_{0}^{t} \Vert e^{(t-\tau)\partial_1^2}e^{\sigma (t-\tau) \mathcal{R}_1 } \partial_1^\alpha  Q^{\theta}_{\sigma}  (\partial_2\mathbb{P}(u \cdot \nabla u), \partial_2(u\cdot \nabla \theta)) \Vert_{L^{2}} d\tau \\
	\lesssim  &  \Vert e^{t\partial_1^2 } \partial_1^\alpha  (\partial_2u_0 | \partial_2\theta_0) \Vert_{L^{2}}+ \int_{0}^{t} \Vert e^{(t-\tau)\partial_1^2 } \partial_1^\alpha  \partial_2(u\cdot \nabla \theta) \Vert_{L^{2}} d\tau  \\
	& + \int_{0}^{t} \Vert e^{(t-\tau)\partial_1^2 }  \partial_1^\alpha  \partial_2(u \cdot \nabla u_1) \Vert_{L^{2}} d\tau  + \int_{0}^{t} \Vert e^{(t-\tau)\partial_1^2 } \partial_1^\alpha  \partial_2(u \cdot \nabla u_2) \Vert_{L^{2}} d\tau  \\
	\triangleq &  \Vert e^{t\partial_1^2 } \partial_1^\alpha  (\partial_2u_0 | \partial_2\theta_0) \Vert_{L^{2}} + H_1 + H_2 + H_3.
	\end{align*}
	The linear term admits a straightforward decay estimate:
	\begin{align*}
	\Vert e^{t\partial_1^2 } \partial_1^\alpha (\partial_2u_0 | \partial_2\theta_0)  \Vert_{L^{2}} & \lesssim (1+t)^{-\frac{1+2\alpha}{4} }(\Vert (\partial_2u_0 | \partial_2\theta_0) \Vert_{ L^{2}_{x_2} (L^{1}_{x_1}) } + \Vert \partial_1^\alpha  (\partial_2u_0 | \partial_2\theta_0) \Vert_{L^2} ) \\
	& \lesssim (1+t)^{-\frac{1+2\alpha}{4}}(\Vert (u_0 | \theta_0) \Vert_{W^{2,1}} + \Vert (u_0 | \theta_0) \Vert_{H^2} ) \\
	& \lesssim  \varepsilon  (1+t)^{-\frac{1+2\alpha}{4} }.
	\end{align*}
	For the nonlinear term $H_1$, by Lemma \ref{HEAT DECAY LEMMA} and the divergence-free condition of $u$, we have 
	\begin{align*}
H_1 \lesssim & \int_{0}^{t} \Vert e^{(t-\tau)\partial_1^2 }  \partial_1^\alpha  \partial_2(u_1 \partial_1 \theta+u_2 \partial_2  \theta) \Vert_{L^{2}} d\tau \\
	\lesssim & \int_{0}^{t}(t - \tau )^{-\frac{1+4\alpha}{8}-\eta}( \Vert   u_1 \partial_1 \partial_2 \theta  \Vert_{ L^2_{x_2} (L^{ \frac{4}{3+8 \eta}}_{x_1})} +\Vert   u_2 \partial_2^2 \theta  \Vert_{ L^2_{x_2} (L^{ \frac{4}{3+8 \eta}}_{x_1})} ) d\tau \\
	& +  \int_{0}^{t}(t - \tau )^{-\frac{1+4\alpha}{8}-\eta}( \Vert  \partial_2 u_1 \partial_1 \theta  \Vert_{ L^2_{x_2} (L^{ \frac{4}{3+8 \eta}}_{x_1})} +\Vert   \partial_1 u_1 \partial_2 \theta  \Vert_{ L^2_{x_2}( L^{ \frac{4}{3+8 \eta}}_{x_1})} ) d\tau \\
	\triangleq & \int_{0}^{t}(t - \tau )^{-\frac{1+4\alpha}{8}-\eta}(H_{11} + H_{12} + H_{13} + H_{14})(\tau) d\tau.
	\end{align*}
	By  Minkowski's inequality, H\"older's inequality and interpolation,
	\begin{align*}
H_{11} (\tau)	\lesssim & \Big\Vert \Vert u_1 \Vert_{L^{4}_{x_2}} \Vert \partial_1 \partial_2 \theta \Vert_{L^{4}_{x_2}} \Big\Vert_{L^{ \frac{4}{3+8 \eta}}_{x_1}}  \\
	\lesssim &  \Vert u_1 \Vert_{L^{\frac{4}{1+8 \eta}}_{x_1} (L^{4}_{x_2})} \Vert \partial_1 \partial_2\theta \Vert_{L^{2}_{x_1}(L^{4}_{x_2}) }  \\
	\lesssim & \Vert u_1 \Vert_{L^{\frac{4}{1+8 \eta}}}^{1-2\eta} \Vert \partial_2 u_1 \Vert_{L^{\frac{4}{1+8 \eta}}}^{2\eta} \Vert \partial_1 \partial_2\theta \Vert_{L^{2}}^{\frac{3}{4}} \Vert \partial_1 \partial_2^2 \theta \Vert_{L^{2}}^{\frac{1}{4} } \\
	\lesssim & \Vert u_1 \Vert_{L^{4}}^{(1-2\eta)(1-8\eta)}\Vert u_1 \Vert_{L^{2}}^{(1-2\eta)8\eta} \Vert \partial_2 u_1 \Vert_{L^{4}}^{2\eta(1-8\eta)} \\
	& \times \Vert \partial_2 u_1 \Vert_{L^{2}}^{2\eta \cdot 8\eta}\Vert \partial_1 \partial_2\theta \Vert_{L^{2}}^{\frac{3}{4}} \Vert \partial_1 \partial_2^2 \theta \Vert_{L^{2}}^{\frac{1}{4} }\\
	\lesssim & \Vert u_1 \Vert_{L^{4}}^{(1-2\eta)(1-8\eta)}\Vert u_1 \Vert_{L^{2}}^{(1-2\eta)8\eta} \Vert  u_1 \Vert_{L^{4}}^{2\eta(1-8\eta)\frac{k-2}{k-1}} \Vert  \partial_2^{k-1} u_1 \Vert_{L^{4}}^{2\eta(1-8\eta)\frac{1}{k-1}}\\
	& \times \Vert \partial_2 u_1 \Vert_{L^{2}}^{2\eta \cdot 8\eta}\Vert \partial_1 \partial_2\theta \Vert_{L^{2}}^{\frac{3}{4}} \Vert \partial_1 \partial_2 \theta \Vert_{L^{2}}^{\frac{k-3}{4k-8} } \Vert \partial_1 \partial_2^{k-1} \theta \Vert_{L^{2}}^{\frac{1}{4k-8} }\\
	\lesssim & C_0^2 \varepsilon^2(1+\tau)^{-m_1},
	\end{align*}
	where $$m_1 = [\frac{1}{2}(1-8\eta) + (\frac{1}{8} + \eta)8\eta](1-2\eta) + [\frac{1}{2}(1-8\eta)(\frac{k-2}{k-1}) + (\frac{1}{8}+ \eta) 8\eta]2\eta + (\frac{5}{8}+ \eta)(\frac{3}{4} + \frac{k-3}{4k-8}).$$ 
	For $\eta = \frac{1}{120}$, $k=14$ and $\delta = 10^{-5}$, we can get  $m_1 \thickapprox 1.0950$  which is larger than 1.\\
	Based on Minkowski's inequality, H\"older's inequality and related interpolation tools, we get
	\begin{align*}
	H_{12} (\tau) 
	\lesssim & \Big\Vert \Vert u_2 \Vert_{L^{\infty}_{x_2}} \Vert \partial_2^2 \theta \Vert_{L^{2}_{x_2}} \Big\Vert_{L^{\frac{4}{3+8 \eta}}_{x_1}} \\
	\lesssim &  \Vert u_2 \Vert_{L^{\frac{4}{1+8\eta}}_{x_1}(L^{\infty}_{x_2} )} \Vert \partial_2^2 \theta \Vert_{L^2}  \\
	\lesssim &  \Vert u_2 \Vert_{L^{\frac{4}{1+8\eta}}_{x_1}(L^{\infty}_{x_2} )} \Vert \partial_2^2 \theta \Vert_{L^2}  \\
	\lesssim &   \Vert u_2 \Vert_{L^{4}}^{(\frac{3}{4} - 2 \eta) (1-8\eta)}  \Vert u_2 \Vert_{L^{2}}^{(\frac{3}{4}- 2 \eta ) 8 \eta} \Vert \partial_1 u_1 \Vert_{L^{4}}^{(\frac{1}{4}+ 2\eta) (1-8\eta)}   \Vert \partial_1 u_1 \Vert_{L^{2}}^{(\frac{1}{4} + 2\eta )  8\eta}\Vert  \partial_2^2  \theta \Vert_{L^{2}}  \\
	\lesssim &   \Vert u_2 \Vert_{L^{4}}^{(\frac{3}{4} - 2 \eta) (1-8\eta)}  \Vert u_2 \Vert_{L^{2}}^{(\frac{3}{4}- 2 \eta ) 8 \eta} \Vert \partial_1 u_1 \Vert_{L^{4}}^{(\frac{1}{4}+ 2\eta) (1-8\eta)}   \Vert \partial_1 u_1 \Vert_{L^{2}}^{(\frac{1}{4} + 2\eta )  8\eta}\Vert  \partial_2  \theta \Vert_{L^{2}}^{\frac{k-2}{k-1}}  \Vert  \partial_2^k  \theta \Vert_{L^{2}}^{\frac{1}{k-1}} \\
	\lesssim &  C_0^2 \varepsilon^2 (1+\tau)^{-m_2},
	\end{align*}
	where $$m_2 = [(\frac{7}{8} + \eta - \delta)(1-8\eta)+ 2\eta](\frac{3}{4} - 2\eta) + [(1-8\eta)(1-\delta)+ 5\eta + 8\eta^2](\frac{1}{4} + 2\eta) + \frac{k-2}{k-1}(\frac{1}{8}+ \eta).$$ 
	For $\eta = \frac{1}{120}$, $k=14$ and $\delta = 10^{-5}$, we can get  $m_2 \thickapprox 1.000030553$  which is larger than 1.
In view of Minkowski's inequality, H\"older's inequality and interpolation,
	\begin{align*}
	H_{13}(\tau)
	\lesssim & \Big\Vert \Vert \partial_1\theta \Vert_{L^{4}_{x_2}} \Vert \partial_2 u_1 \Vert_{L^{4}_{x_2}} \Big\Vert_{L^{ \frac{4}{3+8 \eta}}_{x_1}}  \\
	\lesssim &  \Vert \partial_1 \theta \Vert_{L^{\frac{4}{1+8 \eta}}_{x_1}(L^{4}_{x_2})} \Vert \partial_2 u_1 \Vert_{L^{2}_{x_1}(L^{4}_{x_2} )}  \\
	\lesssim & \Vert \partial_1\theta \Vert_{L^{\frac{4}{1+8 \eta}}}^{1-2\eta} \Vert \partial_1\partial_2 \theta \Vert_{L^{\frac{4}{1+8 \eta}}}^{2\eta} \Vert \partial_2u_1 \Vert_{L^{2}}^{\frac{3}{4}} \Vert  \partial_2^2 u_1 \Vert_{L^{2}}^{\frac{1}{4} } \\
	\lesssim & \Vert \partial_1\theta \Vert_{L^{4}}^{(1-8\eta)(1-2\eta)}\Vert \partial_1\theta \Vert_{L^{2}}^{(1-2\eta)8\eta} \\
	& \times \Vert  \partial_1 \partial_2\theta \Vert_{L^{4}}^{2\eta(1-8\eta)}\Vert  \partial_1 \partial_2 \theta\Vert_{L^{2}}^{8\eta \cdot 2\eta}\Vert \partial_2u_1  \Vert_{L^{2}}^{\frac{3}{4}} \Vert \partial_2^2u_1  \Vert_{L^{2}}^{\frac{1}{4} }\\
	\lesssim & \Vert \partial_1\theta \Vert_{L^{4}}^{(1-8\eta)(1-2\eta)}\Vert \partial_1\theta \Vert_{L^{2}}^{(1-2\eta)8\eta}\Vert  \partial_1 \theta \Vert_{L^{4}}^{2\eta(1-8\eta)\frac{k-3}{k-2}} \Vert  \partial_1 \partial_2^{k-2} \theta \Vert_{L^{4}}^{2\eta(1-8\eta)\frac{1}{k-2}} \\
	& \times \Vert  \partial_1 \partial_2 \theta\Vert_{L^{2}}^{8\eta \cdot 2\eta}\Vert \partial_2u_1  \Vert_{L^{2}}^{\frac{3}{4}} \Vert \partial_2u_1  \Vert_{L^{2}}^{\frac{k-2}{4k-4} } \Vert \partial_2^ku_1  \Vert_{L^{2}}^{\frac{1}{4k-4} }\\
	\lesssim & C_0^2 \varepsilon^2(1+\tau)^{-m_3},
	\end{align*}
	where $$m_3 = (\frac{5}{8}+ \eta)8\eta + (1-\delta)(1-8\eta) [(1-2\eta) +\frac{k-3}{k-2}2\eta] + (\frac{1}{8}+ \eta)[\frac{3}{4} + \frac{k-2}{4k-4}].$$
	 For $\eta = \frac{1}{120}$, $k=14$ and $\delta = 10^{-5}$, we can get  $m_3 \thickapprox 1.1050$  which is larger than 1.\\
The term $H_{14}$ can be estimated similarly to $H_{13}$, and we obtain that
\begin{align*}
H_{14}(\tau)
\lesssim & C_0^2 \varepsilon^2(1+\tau)^{-m_3}.
\end{align*}
	Therefore by Lemma \ref{decay lem}, we get that
	\begin{equation}\label{H_1}
		\begin{aligned}
H_1 =	& \int_{0}^{t} \Vert e^{\partial_1^2 (t-\tau)}  \partial_2(u \cdot \nabla \theta) \Vert_{L^{2}} d\tau  \\
	\lesssim &  C_0^2 \varepsilon^2 (1+t)^{-\frac{1+4\alpha}{8}-\eta}.
	\end{aligned}
	\end{equation} 
	The  nonlinear terms $H_2$ can be estimated similarly,
	\begin{align}\label{H_2}
	H_2 \lesssim C_0^2 \varepsilon^2 (1+t)^{-\frac{1+4\alpha}{8}-\eta}. 
	\end{align}
	For the nonlinear term $H_3$, since the decay of $u_2$ is faster than that of $u_1$ in both the $L^2$ and $L^4$ norms, the estimate of $H_3$ is easier than that of $H_1$; we omit the details here.
	\begin{align}\label{H_3}
 H_3 \lesssim  C_0^2 \varepsilon^2 (1+t)^{-\frac{1+4\alpha}{8}-\eta} .
	\end{align}
	Combining the estimates \eqref{H_1}, \eqref{H_2}, \eqref{H_3} above, we obtain
	\begin{align*}
	\Vert \partial_1^\alpha (\partial_2 u_1 | \partial_2\theta) \Vert_{L^{2}} \lesssim(1+C_0^2 \varepsilon) \varepsilon (1+t)^{-\frac{1+4\alpha}{8} - \eta}.
	\end{align*}
Taking $\varepsilon$ sufficiently small, one can derive the bootstrap
conclusion (\ref{IDpartial2u}).

\vskip .2in
\section*{Acknowledgements}
J. Wu was partially supported by the National Science Foundation of USA (DMS 2104682 and DMS 2309748). N. Zhu was partially supported by the National Natural Science Foundation of China (Grant No. 12301285 and No. 12171010) and by the Shandong Provincial Natural Science Foundation (Project ZR2023QA002).

\vskip .1in
\noindent\textbf{Data availability statement.} Data sharing is not applicable to this article, as no datasets were generated or analyzed during the current study.

\vskip .1in
\noindent\textbf{Conflict of interest.} The authors declare that they have no conflict of interest.

\end{spacing}
\end{document}